\documentclass[a4paper,11pt,spanish, twoside, leqno]{tfg-uam}

\usepackage[utf8]{inputenc}
\usepackage{amsfonts, amssymb, amsmath, amsthm}
\usepackage{graphicx}
\usepackage{color}
\usepackage{verbatim} 
\usepackage{enumitem}
\usepackage{bbm}
 
\newcommand{\rojo}[1]{\textcolor{red}{#1}}
\newcommand*{\avint}{\mathop{\ooalign{$\int$\cr$-$}}}

\newtheorem{teor}{Teorema}[chapter]
\newtheorem{lema}[teor]{Lema}
\newtheorem{prop}[teor]{Proposición}

\newtheorem{prob}[teor]{Problema}
\newtheorem{supo}[teor]{Suposición}

\newtheorem{cor}[teor]{Corolario}

\theoremstyle{definition}

\include{thebibliography}

\title{La ecuación de Keller-Segel }
\author{Alejandro Fernández Jiménez }
\tutor{María del Mar González Nogueras}
\curso{2020-2021}

\usepackage{hyperref}
\hypersetup{
	pdfinfo={
            Title={Las ecuaciones de Keller-Segel \rojo{(Título provisional} }, 
            Author={Alejandro Fernández Jiménez }, 
            Director1={mariamar.gonzalezn }, 
            Director2={Ninguno }, 
            Ndirectores={1 }, 
            Tipo={TFG}, 
            Curso={2018-19}, 
            Palabrasclave={Keller Segel },
				}
}

\begin{document}

\begin{abstract}[spanish]
El propósito de estas trabajo se focaliza en el estudio de la \textit{quimiotaxis} y cómo modelizarlo a través de las ecuaciones de Keller-Segel. La \textit{quimiotaxis} es un proceso natural que provoca que los organismos dirijan sus movimientos atendiendo a ciertos químicos concentrados en su entorno. Es un fenómeno observable para algunas myxamoebas, que experimentan un paseo aleatorio para esparcirse por el espacio buscando comida mientras se ven afectadas por la \textit{quimiotaxis}. El resultado de combinar ambos comportamientos da lugar a una competición entre difusión y agregación. Estas notas pretenden recopilar los resultados más relevantes sobre estas ecuaciones proporcionando una bibliografía exhaustiva sobre los mismos. De esta forma, el primer objetivo es presentar el modelo en EDPs para la densidad de las partículas, planteado por E.F. Keller y L.A. Segel y por C.S. Patlak. Después, hacemos una deducción probabilística de la ecuación. Tras esto, nos damos cuenta de que la masa permanece constante a lo largo del tiempo. En dimensión $2$, y tras simplificar ligeramente el modelo, observamos varios comportamientos según la masa inicial $M$: Para $M< 8 \pi$, domina la difusión, para $M >8 \pi$ (y segundo momento inicial finito) lo hace la agregación, y si $M=8 \pi$ (y segundo momento inicial finito), se produce blow-up en tiempo infinito. Por último, motivado por el empleo de técnicas de entropía en varias partes del trabajo, decidimos desarrollarlas un poco más en otros tipos de ecuaciones.
\end{abstract}
\begin{abstract}[english]
The purpose of this work is the study of \textit{chemotaxis} and how to model it through the equations of Keller-Segel. \textit{Chemotaxis} is a natural process which induces the organisms to direct their movement according to certain chemicals concentrated on their surroundings. This was observed in some myxamoebas, which experiment a random walk, spreading through space, looking for food, while they are affected by chemotaxis. As a result of combining both behaviours, one obtains a competition between diffusion and aggregation. This notes try to make a summary of the main results about these equations, providing an exhaustive bibliography about them. First of all, we would like to present the model in PDEs, we will use the model set by E.F. Keller and L.A. Segel and by C.S. Patlak. Following it, we will do a probabilistic deduction in order to check that these equations could model this situation. After that, we realise that the mass remains constant along time. In dimension $2$, on a slightly simplified model, we observe several behaviours according to the initial mass $M$: If $M<8 \pi$, the diffusion dominates, if $M>8 \pi$ (and finite initial second moment) there will be aggregation, and if $M=8 \pi$ (and finite initial second moment) we get a blow-up in infinite time. Finally, motivated by the use of entropy techniques on several parts of the work, we decide to further develop them on other type of equations.
\end{abstract}

\mainmatter

\chapter{Introducción}
\setcounter{page}{1}

El sistema de Keller-Segel se puede entender como un primer paso para comprender como, durante la evolución de las especies, organismos unicelulares evolucionaron hasta estructuras más complejas. Sirve también como modelo para explicar los patrones de formación de las células por meiosis \cite{CAL_HAW_MEU_VOI}, embrio-genésis o angio-genésis \cite{CHAP, LEV_NIL_SLE}, la enfermedad de Balo \cite{KHO_CAL} y bio-convección \cite{CHE_FEL_KUR_LOR_MAR}. En física este sistema sirve, por ejemplo, para modelizar el movimiento del campo medio de varias partículas Brownianas autogravitatorias, \cite{CHA, CHA_MAN}.

\begin{figure}[h]
  \centering
   \includegraphics[width=1\textwidth]{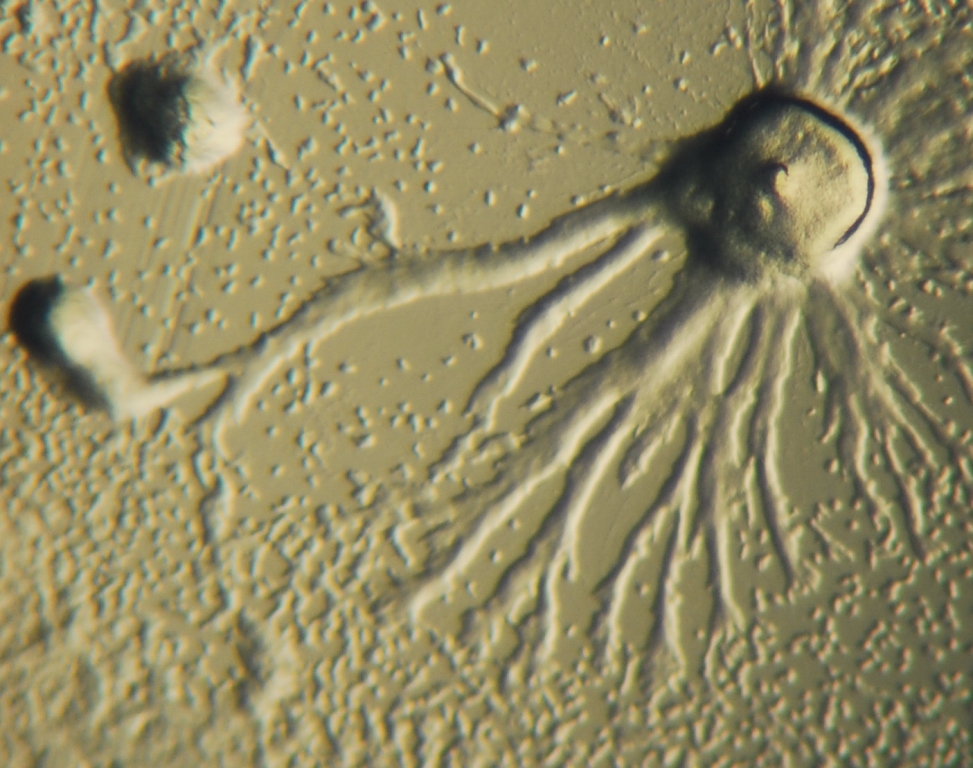}
  \caption{Fenómeno de agregación observado en \textit{Dictostelium discoideum} para formar una estructura más compleja. Fuente: \textit{Wikipedia}.}
  \label{fig:Dictyostelium_Aggregation}
\end{figure}

Este trabajo se centra en el fenómeno de la \textit{quimiotaxis} y cómo modelizarlo a través de las ecuaciones de Keller-Segel. La \textit{quimiotaxis} es un fenómeno que provoca que los organismos dirijan sus movimientos atendiendo a ciertos químicos concentrados en su entorno, la película \url{https://n9.cl/vntg2}\footnote{Investigación sobre la Dictyostelium discoideum del Prof. John Bonner, visualizar 0:14-0:47.} a la Dictyostelium discoideum, un tipo de moho mucilaginoso, ayuda a entender y visualizar estos conceptos. Véase también la figura \ref{fig:Dictyostelium_Aggregation}, en la que se puede observar una fotografía en la que los \textit{Dictostelium discoideum} se han agregado para formar una estructura más compleja. Si el movimiento se dirige hacia concentraciones más altas de la sustancia química hablamos de \textit{quimiotaxis} positiva y el atractor se denomina \textit{quimioatractor}. La evidencia biológica de la existencia de \textit{quimiotaxis} durante la agregación myxobacterial es discutida por Dworkin y Kaiser \cite{DWO_KAI}. 

Algunas células generan ellas mismas el \textit{quimioatractor}, creando interacciones de largo rango no-locales entre ellas. Estas notas pretenden estudiar un modelo de agregación simplificado, entendiendo cómo afecta la \textit{quimiotaxis} a escala celular. Esto es observable en algunas myxamoebas (como por ejemplo las que mencionamos arriba), que experimentan un paseo aleatorio para esparcirse por el espacio y buscar comida mientras se ven afectadas por la \textit{quimiotaxis}. El resultado de combinar ambos fenómenos da lugar a una competición entre difusión y agregación. La difusión es de tipo medios porosos, que ya ha sido ampliamente estudiada por, entre otros, Vázquez (e.g. \cite{VAZ, VAZ2}). 

Los primeros en presentar un modelo matemático para explicar el fenómeno de la \textit{quimiotaxis} fueron E. F. Keller y L. A. Segel en \cite{KEL_SEG} y C. S. Patlak en \cite{PAT}. El resultado fue la siguiente ecuación de drift-difusión que ahora se conoce como el sistema parabólico de Keller-Segel (KS)
\begin{equation} \label{ec:Keller_Segel_Primer_sistema}
\left\{ \begin{array}{lll}
\dfrac{\partial \rho}{\partial t} = \Delta \left( \rho^m \right) - \nabla \cdot \left( \rho \nabla c \right), \\
\varepsilon \partial_t c = \Delta c - \alpha c + \rho, & \, \, \, \, \, \, \, \, \, \, \, \, \, \, \left( t, x \right) \in \left( 0, + \infty \right) \times \mathbb{R}^d, \\
\rho_0 \geq 0, \, \, \, \, \, \, \, c_0 \geq 0.
\end{array} \right.
\end{equation}
En este sistema, $m \in \left[ 1, 2 \right)$, $\varepsilon$ y $\alpha$ son dos parámetros dados no negativos y $d \geq 1$ es la dimensión. $\rho$ representa la densidad celular y $c$ la concentración del quimioatractor ($\rho_0$ y $c_0$ son los valores iniciales). 

Se puede introducir un nuevo modelo que generaliza al anterior,
\begin{equation}\label{ec:difusion_agregacion_caso_general}
\dfrac{\partial \rho}{\partial t} = \Delta \left( \rho^m \right) - \nabla \cdot \left( \rho \nabla \left( \mathcal{K} \ast \rho \right) \right) \, \, \, \, \, \, \mathrm{en} \, \left( 0, + \infty \right) \times \mathbb{R}^d 
\end{equation}
donde $\mathcal{K}$ es un potencial de interacción. El primer término se corresponde con el de difusión y el segundo con el de agregación. Si tomamos $\mathcal{K}$ el núcleo del operador $\varepsilon \partial_t - \Delta + \alpha$ estamos en el caso particular \eqref{ec:Keller_Segel_Primer_sistema}.

Para los Capítulos \ref{ch:Subcritica}, \ref{ch:Critica} y \ref{ch:Supercritica} consideraremos la siguiente versión clásica simplificada de KS. Esta simplificación del modelo fue introducida por primera vez por Nanjundiah \cite{NAN}. Para dimensión $d=2$ es,
\begin{equation} \label{ec:KS_clasico}
\left\{ \begin{array}{lll}
\dfrac{\partial \rho}{\partial t} = \Delta \rho  - \nabla \cdot \left( \rho \nabla c \right) & \, \, \, \, x \in \mathbb{R}^2, \, t>0, \\
- \Delta c =  \rho & \, \, \, \, x \in \mathbb{R}^2, \, t>0,  \\
\rho \left( \cdot , t=0 \right) = \rho_0 \geq 0 & \, \, \, \, x \in \mathbb{R}^2.
\end{array} \right.
\end{equation}
Este es un modelo que puede ser visto como el caso límite para el cuál el atractor químico se difunde mucho más rápido de lo que lo hacen las células que lo emiten.

En la segunda ecuación de \eqref{ec:KS_clasico} elegimos $c$ dada por $c = G \ast \rho$, donde $G$ es el núcleo de Poisson definido por $G \left( \left| x \right| \right) := - \frac{1}{2 \pi} \log \left| x \right|$. De esta manera, el sistema de KS \eqref{ec:KS_clasico} puede ser escrito con la ecuación parabólica no-local:
\begin{equation}\label{ec:KS_clasico_una_linea}
\dfrac{\partial \rho}{ \partial t} = \Delta \rho - \nabla \cdot \left( \rho \nabla G \ast \rho \right) \, \, \, \, \mathrm{en} \, \left( 0, + \infty \right) \times \mathbb{R}^2.
\end{equation}
Varios experimentos sugieren que, con la cantidad suficiente de bacterias, estas se agregan, mientras que, para cantidades bajas, estas se dispersan, e.g. \cite{BRE_CON_KAD_SCH_VEN}. Para aplicar estas ideas, observamos una propiedad de gran relevancia de las ecuaciones de Keller-Segel: Si asumimos suficiente decaimiento en el infinito y regularidad de $\rho$, es fácil demostrar que las soluciones del problema \eqref{ec:difusion_agregacion_caso_general} preservan la masa a lo largo del tiempo 
\begin{equation}\label{ec:def_masa}
M := \int_{\mathbb{R}^d} \rho_0 \left( x \right) dx = \int_{\mathbb{R}^d} \rho \left( t, x \right) dx.  
\end{equation}
Como resultado de juntar todos estos factores, sospechamos que la masa jugará un papel importante en la evolución de la ecuación. Para discutir sobre la influencia que tiene la masa en la solución de \eqref{ec:KS_clasico} definimos el segundo momento
\begin{equation*}
M_2 \left( t \right) := \int_{\mathbb{R}^2} \left| x \right|^2 \rho \left( x, t \right) dx .
\end{equation*}
Suponiendo suficiente regularidad y decaimiento en el infinito de la solución $\rho$, es fácil demostrar utilizando las igualdades anteriormente descritas que
\begin{equation*}
\dfrac{d}{dt} M_2 \left( t \right) = \dfrac{d}{dt} \int_{\mathbb{R}^2} \left| x \right|^2 \rho \left( x, t \right) dx = 4M \left( 1 - \dfrac{M}{8 \pi} \right)
\end{equation*}
con $M$ la masa, \eqref{ec:def_masa}. Entonces, si el segundo momento es finito, tenemos 
\begin{equation*}
\int_{\mathbb{R}^2} \left| x \right|^2 \rho \left( x, t \right) dx = \int_{\mathbb{R}^2} \left| x \right|^2 \rho \left( x, 0 \right) dx + 4M \left( 1 - \dfrac{M}{8 \pi} \right) t.
\end{equation*}
Como consecuencia,
\begin{itemize}
\item Si $M < 8 \pi$ (\underline{masa subcrítica})  el segundo momento crece de manera lineal con el tiempo mientras la masa se conserva. Por tanto, $\rho \left( x, t \right)$ se difunde. En el Capítulo \ref{ch:Subcritica} demostramos que, efectivamente, se produce difusión, empleando el método JKO. Seguimos las ideas expuestas por Blanchet, Calvez y Carrillo en \cite{BLA_CAL_CAR}.
\item Si $M = 8 \pi$ (\underline{masa crítica})  el segundo momento de la solución se preserva con el tiempo pero no es claro que comportamiento adoptarán las soluciones. Por ello, siguiendo un resultado reciente de Dávila, del Pino, Dolbeault, Musso y Wei \cite{DAV_PIN_DOL_MUS_WEI}, estudiamos este caso en el Capítulo \ref{ch:Critica}, para concluir que, con segundo momento inicial finito, se produce blow-up en tiempo infinito, y que estos toman la forma de bubbles (burbujas).
\item Si $M > 8 \pi$ (\underline{masa supercrítica}) la solución no puede permanecer suave a partir de cierto tiempo. $\rho \left( x, t \right)$ explota en tiempo finito. En el Capítulo \ref{ch:Supercritica} estudiamos un ejemplo dado por Herrero y Velázquez \cite{HER_VEL} para el que tenemos blow-up, y donde se desarrolla el método de \textit{matching asymptotics}.
\end{itemize}

Resumiendo, nuestros objetivos en este trabajo son comprender las ecuaciones de Keller-Segel, comenzando por la derivación de estas con técnicas probabilísticas, para ello, en el Capítulo \ref{ch:Probabilidad} derivamos formalmente las ecuaciones de Keller-Segel \eqref{ec:Keller_Segel_Primer_sistema} como límite dinámico de sistemas estocásticos multipartículas con interacciones moderadas. Tras esto, en los Capítulos \ref{ch:Subcritica}, \ref{ch:Critica} y \ref{ch:Supercritica}, intentaremos entender las soluciones de la ecuación para datos iniciales distintos, separando por la masa inicial, como hemos discutido arriba. Por último, motivado por el empleo de las técnicas de entropía en varios lugares del trabajo, decidimos dedicar el Capítulo \ref{ch:entropia} a profundizar en ellas para distintas ecuaciones. Comenzamos entendiendo cómo funcionan para las ecuaciones de Fokker-Planck, y, posteriormente, estudiamos el método en ecuaciones algo más generales. De esta forma, estas notas pretenden proporcionar un resumen sobre los resultados más relevantes acerca de las ecuaciones de Keller-Segel, proporcionando bibliografía exhaustiva sobre estas.

Ademas, acompañamos el trabajo con varios Apéndices, que permiten complementar las distintas técnicas utilizadas a lo largo de estas notas. En el Apéndice \ref{ap:Transporte_optimo} se motiva el empleo de la métrica de Wasserstein, y para ello se presenta el problema del transporte óptimo y se comentan algunas de las aplicaciones dadas por Otto al transporte óptimo. En el Apéndice \ref{ap:HWI}, se enuncian y se demuestran las desigualdades HWI para completar la Sección \ref{sec:Car_McC_Vil_Convergencia_exponencial}, ya que desempeñan un papel muy importante en la demostración del resultado principal de la Sección. En el Apéndice \ref{ap:Liouville}, hablamos sobre la ecuación de curvatura Gaussiana, para entender mejor la ecuación de Liouville, un caso particular de la primera, que aparece en el Capítulo \ref{ch:Subcritica}, y demostrar existencia y unicidad de soluciones para la ecuación desde un punto de vista más geométrico. Finalmente, en el Apéndice \ref{ap:auto-similar}, introducimos las singularidades tipo I y tipo II para una ecuación de evolución por medio de ejemplos. Para las desigualdades tipo I comentamos el caso de una superficie sólida evolucionando bajo la acción de difusión superficial en la que los átomos migran a lo largo de la superficie, conducidos por gradientes de potenciales químicos, y también hablamos sobre los \textit{thin films and thin jets}. Y para estudiar las singularidades tipo II, consideramos la ecuación de Burgers. Este último Apéndice lo redactamos para introducir algunos de los conceptos e ideas más importantes que utilizamos durante el Capítulo \ref{ch:Supercritica}.

\chapter{Derivación probabilística}\label{ch:Probabilidad}

En este Capítulo pretendemos obtener una derivación formal de las ecuaciones de Keller-Segel \eqref{ec:Keller_Segel_Primer_sistema} como límite dinámico de sistemas estocásticos multipartículas con interacciones moderadas. La primera derivación rigurosa de este tipo de sistemas fue presentada por Stevens en \cite{STE} recurriendo a algunas ideas que desarrolló en \cite{STE2}. En este Capítulo nosotros seguiremos \cite{STE}.

Aquí, utilizaremos $u$ para hacer referencia a los términos relacionados con las bacterias y $v$ para los relacionados con sustancias químicas. Sea $S \left( N, t \right) := S_u \left( N, t \right) + S_v \left( N, t \right)$ el conjunto de todas las partículas de nuestro sistema, donde $N$ es el número total de partículas del sistema, que son numeradas tomando un nuevo número para cada partícula recién nacida. Definimos $P_N^k \left( t \right) \in \mathbb{R}^d$, $k \in S \left( N, t \right)$ para poder describir la posición de la partícula $k$ en tiempo $t \geq 0$. También definimos la siguiente medida evaluada en procesos empíricos
\begin{equation*}
t \mapsto S_{N_r} ( t ) = \frac{1}{N} \sum_{k \in S_{N_r} ( N, t )} \delta_{P_N^k \left( t \right)} \; \; \; \;  \mathrm{y} \; \; \; \; S_N ( t ) = S_{N_u} ( t ) + S_{N_v} ( t ),
\end{equation*} 
$r=u$, $v$. Las funciones que describen la dinámica de las partículas dependen de las siguientes versiones suavizadas de $S_{N_r} ( t )$
\begin{equation*}
\hat{s}_{N_r} ( t, x ) = ( S_{N_r} ( t ) \ast W_N \ast \hat{W}_N ) ( x )
\end{equation*}
donde $W_N$ y $\hat{W}_N$ son densidades de probabilidad, obtenidas a partir de una función suave simétrica fijada de antemano $W_1$,
\begin{equation*}
W_N \left( x \right) = \alpha_N^d W_1 \left( \alpha_N x \right) \; \; \; \;  \mathrm{y} \; \; \; \; \hat{W}_N \left( x \right) = \hat{\alpha}_N^d W_1 \left( \alpha_N x \right),
\end{equation*}
con $\alpha_N = N^{\alpha/d}$, $\hat{\alpha}_N = N^{\hat{\alpha}/d}$ y $\alpha$, $\hat{\alpha}$ constantes. De esta forma, el movimiento de cada partícula a través del espacio viene descrito por una ecuación diferencial estocástica obtenida de modelizar un grupo de myxobacterias de tamaño $8 \times 1$ unidades, 
\begin{equation*}
\begin{array}{ll}
dP^k_N \left( t \right) = \chi_N \left( t, P^k_N ( t ) \right) \nabla \hat{s}_{N_v} \left( t, P_N^k ( t ) \right) dt + \sqrt{2 \mu} \, dW^k \left( t \right) & \; \; \mathrm{si} \; \; k \in S_u ( N, t ) \\
 dP^k_N \left( t \right) = \sqrt{2 \eta } \, dW^k ( t ) & \; \; \mathrm{si} \; \; k \in S_v ( N, t ) 
\end{array} 
\end{equation*}
donde $W^k \left( \cdot \right)$ son movimientos Brownianos independientes $d$-dimensionales, $\mu$, $\eta > 0$ son constantes y $\chi_N \left( t, x \right) = \chi \left( \hat{s}_{N_u} \left( t, x \right), \hat{s}_{N_v} \left( t, x \right) \right)$, con $\chi$ un funcional que denota la sensibilidad quimiotáctica de las bacterias.

Mirando a las dinámicas de las partículas de las dos subpoblaciones, uno utiliza la fórmula de Itô para calcular $ f \in C_b^{1, \, 2} \left( \mathbb{R}_{+} \times \mathbb{R}^d \right)$, que es solución de
\begin{equation}\label{ec:S_u,f}
\begin{array}{llll}
\left\langle S_{N_u} \left( t \right), f \left( t, \cdot \right) \right\rangle & = \frac{1}{N} \sum_{k \in S_{N_u} \left( N, t \right)} f \left( t, P_N^k \left( t \right) \right) = \left\langle S_{N_u} \left( 0 \right), f \left( 0, \cdot \right) \right\rangle \\
 & + \int_0^t < S_{N_u} \left( \tau \right), \chi_N \left( \tau , \cdot \right) \nabla \hat{s}_{N_v} \left( \tau , \cdot \right) \cdot \nabla f \left( \tau , \cdot \right) \\
 & + \mu \Delta f \left( \tau , \cdot \right) + \frac{\partial}{\partial \tau} f \left( \tau , \cdot \right) > d \tau \\
 & + \frac{1}{N} \int_0^t \sum_{k \in S_{N_u} \left( N, \tau \right)} \sqrt{2 \mu} \, \nabla f \left( \tau , P_N^k \left( \tau \right) \right) \cdot d W^k \left( \tau \right) 
\end{array}
\end{equation}
y 
\begin{equation}\label{ec:S_v,f}
\begin{array}{lll}
\left\langle S_{N_v} \left( t \right) , f \left( t, \cdot \right) \right\rangle & = \frac{1}{N} \sum_{k \in S_{N_v} \left( N, t \right)} f \left( t, P_N^k \left( t \right) \right) = \left\langle S_{N_v} \left( 0 \right) , f \left( 0, \cdot \right) \right\rangle \\
 & + \int_0^t < S_{N_v} \left( \tau \right) , \eta \Delta f \left( \tau , \cdot \right) + \dfrac{\partial}{\partial \tau} f \left( \tau , \cdot \right) > d \tau \\
 & + \frac{1}{N} \int_0^t \sum_{k \in S_{N_v} \left( N, \tau \right)} \sqrt{2 \eta} \, \nabla f \left( \tau , P_N^k \left( \tau \right) \right) \cdot d W^k \left( \tau \right) 
 \end{array}
\end{equation}
Definimos
\begin{equation}\label{ec:Martingalas_Deduccion_Prob}
M_{N_r}^1 \left( t, f \right) = \frac{1}{N} \int_0^t \sum_{k \in S_{N_r} \left( N, \tau \right)} \sqrt{2 \zeta_r} \, \nabla f \left( \tau , P_N^k \left( \tau \right) \right) \cdot dW^k \left( \tau \right) \\
\end{equation}
donde $r=u$, $v$ y $\zeta_u = \mu$, $\zeta_v = \eta$. Estos dos procesos son martingalas con respecto a la filtración natural $ \left\lbrace \mathcal{F}_t \right\rbrace_{t \geq 0 }$ generadas por los procesos $t \rightarrow \left( P_N^k \left( t \right), \mathbbm{1}_{S_r \left( N, t \right)} \right)$, $r=u$, $v$. Por la derivación heurística del comportamiento límite de \eqref{ec:S_u,f} y \eqref{ec:S_v,f}, uno asume para $r =u$, $v$ y $t \geq 0$
\begin{equation*}
\lim_{N \to \infty} S_{N_r} \left( t \right) = S_r \left( t \right) \; \; \; \; \mathrm{en \; algun \; sentido  ,}
\end{equation*}
donde la medida $S_r \left( t \right)$ tiene una densidad suave $r \left( t , \cdot \right)$ con respecto a la medida de Lebesgue. Sean $u_0 \left( \cdot \right)$ y $v_0 \left( \cdot \right)$ las densidades de $S_u \left( 0 \right)$ y $S_v \left( 0 \right)$. Dado que,
\begin{equation*}
\lim_{N \to \infty } W_N = \lim_{N \to \infty} \hat{W}_N = \delta_0 \; \; \; \; \mathrm{en \; el \; sentido \; de \; distribuciones ,}
\end{equation*}
se sigue que $\lim_{N \to \infty} \hat{s} \left( t , \cdot \right) = r \left( t, \cdot \right)$ y $\lim_{N \to \infty} \nabla \hat{s}_{N_r} \left( t, \cdot \right) = \nabla r \left( t, \cdot \right)$. Asumiendo que las variaciones cuadráticas de $M^1_{N_u} \left( \cdot , f \right)$ y $M^1_{N_v} \left( \cdot , f \right)$ tienden a cero para $N \to \infty$, usando \eqref{ec:S_u,f}, \eqref{ec:S_v,f} y \eqref{ec:Martingalas_Deduccion_Prob} se puede obtener de manera formal, que para $N \to \infty$, 
\begin{equation*}
\begin{array}{ll}
\left\langle u \left( t, \cdot \right), f \left( t, \cdot \right) \right\rangle = & \left\langle u_0 \left( \cdot \right) , f \left( 0 , \cdot \right) \right\rangle + \int_0^t \left\langle u \left( \tau , \cdot \right), \chi_{\infty} \left( \tau , \cdot \right) \nabla v \left( \tau , \cdot \right) \cdot \nabla f \left( \tau , \cdot \right) \right\rangle d \tau \\
 & + \int_0^t \left\langle u \left( \tau , \cdot \right) , \mu \Delta f \left( \tau , \cdot \right) + \frac{\partial}{\partial \tau} f \left( \tau , \cdot \right) \right\rangle d \tau \, , \\
\left\langle v \left( t, \cdot \right), f \left( t, \cdot \right) \right\rangle = &  \left\langle v_0 \left( \cdot \right) , f \left( 0 , \cdot \right) \right\rangle + \int_0^t \left\langle v \left( \tau , \cdot \right) , \eta \Delta f \left( \tau , \cdot \right) + \frac{\partial}{\partial \tau} f \left( \tau , \cdot \right) \right\rangle d \tau \, ,
\end{array}
\end{equation*}
donde $\chi_{\infty} \left( \tau , x \right) = \chi \left( u \left( \tau , x \right) , v \left( \tau , x \right) \right)$. Si Integramos por partes según el sentido de la fórmula de Itô conseguimos una versión simplificada del sistema de \textit{quimiotaxis}
\begin{equation}
\left\{
\begin{array}{ll}
u_t & = \nabla \cdot \left( \mu \nabla u - \chi \left( u, v \right) u \nabla v \right) = \mu \Delta u - \nabla \cdot \left( \chi \left( u, v \right) u \nabla v \right) \\
v_t & = \eta \Delta v
\end{array}
\right.
\end{equation}
donde $u \left( 0 , x \right) = u_0 \left( x \right)$ y $v \left( 0, x \right) = v_0 \left( x \right)$. Con los cálculos hechos hasta ahora solo hemos tendríamos soluciones débiles. Para obtener soluciones fuertes se necesitan algunos detalles técnicos extra. De esta forma, hemos deducido la ecuación de Keller-Segel \eqref{ec:Keller_Segel_Primer_sistema} para $m=1$ salvo por la existencia de dos términos lineales en la evolución temporal de la concentración de químicos. Estos se consiguen añadiendo en nuestros argumentos la posibilidad de que nazcan o mueran nuevas partículas. Todos los pasos extra que comentamos en este párrafo los discute Stevens en \cite{STE}.

\chapter{Masa subcrítica}\label{ch:Subcritica}
 
Bajo la hipótesis de masa subcrítica $M \in (0, 8 \pi )$, el segundo momento crece de manera lineal con el tiempo mientras la masa se conserva. Esto sugiere que la solución se difunde y no se producen blow-ups. Comenzamos definiendo la \textit{energía libre} $E_{KS}$.  
\begin{equation*}
E_{KS} \left[ \rho \right] := \int_{\mathbb{R}^2} \rho \log \rho \, dx - \frac{1}{2} \int_{\mathbb{R}^2} \rho c \, dx ,
\end{equation*}
donde la primera parte de $E_{KS}$, $\int_{\mathbb{R}^2} \rho \log \rho \, dx$, es un término de \textit{entropía} y la segunda, $- \frac{1}{2} \int_{\mathbb{R}^2} \rho c \, dx$, es un \textit{potencial}. El funcional $E_{KS}$ pertenece estructuralmente a la clase general de energías libres para partículas en interacción, introducido en \cite{MCC, CAR_MCC_VIL, CAR_MCC_VIL_2} por Carrillo, McCann y Villani. 

Más adelante, veremos que $E_{KS}$ es el flujo gradiente del sistema de Keller-Segel \eqref{ec:KS_clasico}. Aprovechando esto, vamos a construir una sucesión minimizante de $E_{KS}$, para convencernos de que existe convergencia de las soluciones al estado estacionario (que es el minimizante de la sucesión) en tiempo. Todo esto lo haremos siguiendo las notas de Blanchet \cite{BLA}, utilizando el método de Jordan-Kinderlehrer-Otto (JKO) \cite{JOR_KIN_OTT} y el trabajo de McCann \cite{MCC}. Usamos en concreto, las ideas expuestas en \cite{BLA_CAL_CAR} por Blanchet, Calvez y Carrillo. 


Dividimos este Capítulo en cuatro Secciones. En la primera, utilizamos el término de entropía de la energía libre, dándonos cuenta de que no es suficiente para resolver el problema en todo su rango de valores, $M \in \left( 0, 8 \pi \right)$. No obstante, es una buena introducción e indica que necesitaremos desarrollar técnicas para dar una demostración completa, sugiriendo la necesidad de involucrar todos los términos de $E_{KS}$. Por ello, la siguiente Sección estará dedicada a explicar y entender la relación entre la energía libre $E_{KS}$ y el sistema de Keller-Segel \eqref{ec:KS_clasico}. En la tercera Sección explicamos como se puede utilizar y adaptar el método JKO a este problema particular, haciendo analogías con la estructura del flujo gradiente en el caso euclídeo para facilitar la comprensión del método. En la cuarta Sección damos las ideas principales para probar el Teorema \ref{th:Masa_subcritica_JKO}, que permite juntar toda esta información y completar la demostración que queremos sobre convergencia. Además, en \cite{BLA_CAL_CAR}, Blanchet, Calvez y Carrillo consiguen demostrar convergencial exponencial para el caso unidimensional, pero nosotros no lo discutiremos aquí. 

\section{Primera aproximación al problema}\label{sec:Subcritica_Primera_Seccion}
Lo que a primera vista parece más natural para intentar resolver el problema es estudiar el signo de la variación del término de la entropía de $E_{KS}$. Los primeros resultados siguiendo esta línea son gracias a Jäger y Luckhaus en \cite{JAG_LUC}. 

Nosotros seguimos \cite{COR_PER_ZAA, COR_PER_ZAA_2}, donde se comienza calculando $\frac{d}{dt} \int_{\mathbb{R}^2} \rho \log \rho \, dx$. Si además se integra por partes y se aprovecha la expresión para $c$, obtenemos,
\begin{equation*}
\frac{d}{dt} \int_{\mathbb{R}^2} \rho \log \rho \, dx = -4 \int_{\mathbb{R}^2} \left| \nabla \sqrt{\rho} \right|^2 \, dx + \int_{\mathbb{R}^2} \rho^2 \, dx .
\end{equation*} 
El signo depende de la competición entre dos términos, la difusión basada en el término de \textit{disipación de entropía} $\int_{\mathbb{R}^2} \left| \nabla \sqrt{\rho} \right|^2 dx$ y la \textit{producción hiperbólica de entropía}.

Para cuantificar esta competición podemos usar la desigualdad de Gagliardo-Nirenberg-Sobolev para $p=4$ \cite{GAG, NIR}, donde $C_{GNS} = C_{GNS}^{(4)}$ es la mejor constante para este caso. 
\begin{equation*}
\| u \|^2_{L^p \left( \mathbb{R}^2 \right)} \leq C_{GNS}^{(p)} \| \nabla u \|^{2-4/p}_{L^2 \left( \mathbb{R}^2 \right)} \| u \|^{4/p}_{L^2 \left( \mathbb{R}^2 \right)} \, \, \, \, \forall u \in H^1 \left( \mathbb{R}^2 \right), \, \, \forall p \in \left[ 2 , \infty \right) .
\end{equation*}
Tomando $u=\sqrt{\rho}$ comprobamos que la entropía es no-creciente si $M \leq 4 C_{GNS}^{-2}$. El valor explícito de $C_{GNS}$ no es conocido pero podemos calcularlo numéricamente (véase \cite{WEI}). Convenciéndonos así de que la entropía es no creciente si $M \leq 1.862 \cdots \times (4 \pi ) < 8 \pi$. Esta estimación no cubre todo el rango para el que habíamos predicho que la masa $M$ es subcrítica. La demostración de todo el rango predicho, $M \in (0, 8 \pi)$ fue hecha por Campos y Dolbeault en \cite{CAM_DOL} con técnicas de simetrización adaptadas de un trabajo previo de Díaz, Nagai y Rakotoson \cite{DIA_NAG_RAK}. 
Nosotros en cambio, seguiremos la demostración dada por Blanchet, Calvez y Carrillo \cite{BLA_CAL_CAR}, basada en el uso de flujos gradientes. Daremos un esquema de esta demostración en la Sección \ref{sec:Sub_dem} introduciendo para ello varios conceptos en las Secciones \ref{sec:Sub_energia_libre} y \ref{sec:Sub_JKO}, siguiendo \cite{BLA}. 

\section{Relación entre la energía libre y Keller-Segel}\label{sec:Sub_energia_libre}
Pretendemos estudiar en detalle la energía libre $E_{KS}$ y entender por qué es tan relevante para el problema \eqref{ec:KS_clasico}, de dimensión $d=2$. Un simple cálculo formal demuestra que para todo $u \in C_{c}^{\infty} \left( \mathbb{R}^2 \right)$ con media cero,
\begin{equation*}
\lim_{\epsilon\rightarrow 0} \dfrac{E_{KS} \left[ \rho + \epsilon u \right] - E_{KS} \left[ \rho \right]}{\epsilon} = \int_{ \mathbb{R}^2} \dfrac{\delta E_{KS} \left[ \rho \left( t \right) \right]}{\delta \rho} \left( x \right) u \left( x \right) dx
\end{equation*}
donde
\begin{equation*}
\dfrac{\delta E_{KS} \left[ \rho \right]}{\delta \rho} (x) := \log \rho (x) - G \ast \rho (x) ,
\end{equation*}
con $G$ el potencial de \eqref{ec:KS_clasico_una_linea}. Recordamos que $c=G \ast \rho$ y que $G$ es el núcleo de Poisson, que, para el caso de dimensión $d=2$ toma la forma explícita $G(|x|) = - \frac{1}{2 \pi} \log |x|$. De esta manera, es fácil ver que el sistema \eqref{ec:KS_clasico} puede reescribirse como
\begin{equation*}
\dfrac{\partial \rho}{\partial t} \left( t, x \right) = \nabla \cdot \left( \rho \left( t, x \right) \nabla \left[ \dfrac{\delta E_{KS} \left[ \rho \left( t \right) \right]}{\delta \rho} \left( x \right) \right] \right),
\end{equation*}
de donde se sigue que a lo largo de las trayectorias del sistema KS \eqref{ec:KS_clasico} se satisface
\begin{equation*}
\dfrac{d}{dt} E_{KS} \left[ \rho (t) \right] = - \int_{\mathbb{R}^2} \rho \left( x, t \right) \left| \nabla \left[ \dfrac{\delta E_{KS} \left[ \rho (t) \right]}{\delta \rho} (x) \right] \right|^2 dx.
\end{equation*}
O, lo que es equivalente,
\begin{equation*}
\dfrac{d}{dt} E_{KS} \left[ \rho \left( t \right) \right] = - \int_{\mathbb{R}^2} \rho \left( t, x \right) \left| \nabla \left( \log \rho \left( t, x \right) - c \left( t, x \right) \right) \right|^2 dx.
\end{equation*}
En particular, a lo largo de las trayectorias, $t \mapsto E_{KS} \left[ \rho \left( t \right) \right]$ es monótona no-creciente.

Tras obtener esta relación queremos estudiar su acotación, tal y como hacíamos con nuestro primer método más ingenuo de la Sección \ref{sec:Subcritica_Primera_Seccion}. Esta cota se puede conseguir gracias a la conexión entre la energía libre $E_{KS}$ y la desigualdad logarítmica de Hardy-Littlewood-Sobolev (LogHLS a partir de ahora) que fue conectada por primera vez con esta teoría por Dolbeault y Perthame en \cite{DOL_PER} y que se enuncia:

\begin{prop}
[Desigualdad logarítmica de Hardy-Littlewood-Sobolev, \cite{BEC, CAR_LOS}]  
Sea $f$ una función no-negativa en $\mathcal{L}^1 \left( \mathbb{R}^2 \right)$ tal que $f \log f$, $f \log ( 1+ \left| x \right|^2 ) \in \mathcal{L}^1 \left( \mathbb{R}^2 \right)$ y tal que $\int_{\mathbb{R}^2} f \, dx = M$, entonces
\begin{equation}\label{ec:LogHLS_subcrit}
\int_{\mathbb{R}^2} f \log f \, dx + \dfrac{2}{M} \int \int_{\mathbb{R}^2 \times \mathbb{R}^2 } f \left( x \right) f \left( y \right) \log \left| x- y \right| \, dx \, dy \geq - C \left( M \right),
\end{equation}
con $C \left( M \right) := M \left( 1+ \log \pi - \log M \right)$. Además el minimizador de LogHLS en la desigualdad \eqref{ec:LogHLS_subcrit} son las traslaciones de
\begin{equation*}
f_{\lambda} \left( x \right) := \frac{M}{\pi} \frac{\lambda}{( \lambda + \left| x \right|^2 )^2}.
\end{equation*}
\end{prop}

Usando la monotonía de la energía libre $E_{KS}$ y la desigualdad LogHLS \eqref{ec:LogHLS_subcrit} es fácil ver que para soluciones suaves de \eqref{ec:KS_clasico} tenemos
\begin{equation} \label{ec:Subcritico_Cota_Energia_Libre}
\begin{array}{lll}
E_{KS} \left[ \rho \right] & = \dfrac{M}{8 \pi} \left( {\displaystyle \int}_{\mathbb{R}^2} \rho (x) \log \rho (x) \, dx + \dfrac{2}{M} {\displaystyle \int} {\displaystyle \int}_{\mathbb{R}^2 \times \mathbb{R}^2} \rho (x) \log |x-y| \, dx \, dy  \right) \\
 & \, \, \, \, \, \, \, \, \, \, \, \, \, + \left( 1 - \dfrac{M}{8 \pi} \right) {\displaystyle \int}_{\mathbb{R}^2} \rho (x) \log \rho (x) \, dx \\
 & \geq - \dfrac{M}{8 \pi} C \left(M \right) + \left( 1 - \dfrac{M}{8 \pi} \right) {\displaystyle \int}_{\mathbb{R}^2} \rho \left( x \right) \log \rho \left( x \right) \, dx.
\end{array}
\end{equation}
Se sigue que
\begin{equation} \label{ec:Subcritico_Estimacion}
\int_{\mathbb{R}^2} \rho \left( t, x \right) \log \rho \left( t, x \right) \, dx \leq \dfrac{8 \pi E_{KS} \left[ \rho_0 \right] - M \, C( M)}{8 \pi - M}.
\end{equation}
Entonces, para $M < 8 \pi$, la entropía queda acotada uniformemente en tiempo. Esto, formalmente, impide el colapso de la masa en un punto para un dato y constituirá un argumento crucial en esta demostración.

Merece la pena detenerse un momento para señalar que dado un $\rho$, si fijamos la familia $\rho_{\lambda} (x) = \lambda^{-2} \rho \left( \lambda^{-1} x \right)$, entonces
\begin{equation}
E_{KS} \left[ \rho_{\lambda} \right] = E_{KS} \left[ \rho \right] - 2 M \left( 1 - \dfrac{M}{8 \pi} \right) \log \lambda .
\end{equation}
Así, como función de $\lambda$, la energía $E_{KS} \left[ \rho_{\lambda} \right]$ es acotada por abajo si $M < 8 \pi$, pero no lo es si $M > 8 \pi$ en el conjunto
\begin{equation}\label{set:subcr_funciones_admisibles}
\mathcal{A} := \left\lbrace \rho \, : \, \int_{\mathbb{R}^2} \rho = M, \, \int_{\mathbb{R}^2} \rho (x) \log \rho (x) \, dx < \infty \, \, \mathrm{y} \, \int_{\mathbb{R}^2} | x |^2 \rho (x) \, dx < \infty \right\rbrace  \, .
\end{equation}

\section{Método JKO para el problema de Keller-Segel}\label{sec:Sub_JKO}

Aquí presentamos el método JKO para construir una solución de \eqref{ec:KS_clasico}. No estudiamos Keller-Segel directamente en la forma \eqref{ec:KS_clasico}, sino que lo intetamos reescribir como un flujo gradiente. En vez de hacerlo con la métrica euclídea, nos damos cuenta de que usando la métrica $2$-Wasserstein se tiene,
\begin{equation}\label{ec:flujo_Wasserstein}
\dfrac{\partial \rho}{\partial t} = -  \nabla_{\mathcal{W}_2}  E_{KS} \left[ \rho (t) \right].
\end{equation}
La idea de utilizar esta métrica surge como una de las aplicaciones del transporte óptimo dadas por Otto en \cite{OTT}. Para el lector interesado, puede acudir al Apéndice \ref{ap:Transporte_optimo}, donde encontrará una aproximación más profunda y rigurosa al concepto flujo gradiente (gradient flow) y a la definición de la distancia de Wasserstein $\mathcal{W}_2$. 

Para construir el estado estacionario de \eqref{ec:flujo_Wasserstein} vamos a usar un método de minimización conocido como el método de Jordan-Kinderlehrer-Otto (JKO) \cite{JOR_KIN_OTT}. En este, dado un paso en tiempo $\tau$ definimos la solución por 
\begin{equation}\label{ec:JKO}
\rho_{\tau}^{k+1} \in {\arg\min}_{ \rho \in \mathcal{A}} \left[ \dfrac{\mathcal{W}_2^2 \left( \rho , \rho_{\tau}^k \right)}{2 \tau} + E_{KS} \left[ \rho \right] \right] \, ,
\end{equation}
donde $\mathcal{A}$ está definido en \eqref{set:subcr_funciones_admisibles} y $\mathcal{W}_2$ es la distancia 2-Wasserstein.

Para dar intuición al método JKO, profundizamos en la analogía entre el flujo gradiente para la distancia $2$-Wasserstein y la estructura de este en un caso euclídeo \cite{SAN}. Para este caso, la ecuación de Euler-Lagrange asociada a
\begin{equation}
X_{\tau}^{k+1} \in \arg \min \left[ \dfrac{\left| X - X_{\tau}^k \right|^2}{2 \tau} + E_{KS} \left[ X \right] \right] \, ,
\end{equation}
sería
\begin{equation*}
\dfrac{X_{\tau}^{k+1} - X_{\tau}^k}{ \tau} + \nabla E_{KS} \left[ X_{\tau}^{k+1} \right] = 0 \, ,
\end{equation*}
que no es otra cosa más que el método implícito de Euler para 
\begin{equation*}
\dot{X} = - \nabla E_{KS} \left[ X (t) \right] .
\end{equation*}
Volviendo a nuestro problema, pretendemos construir una sucesión $\left\lbrace \rho_{\tau}^k \right\rbrace_k$ usando la regla \eqref{ec:JKO} y obteniendo en el límite un flujo gradiente que formalmente se puede escribir como en \eqref{ec:flujo_Wasserstein}.

En una disposición euclídea el siguiente paso clásico sería realizar una interpolación entre los nodos. Aquí interpolamos entre los términos de la sucesión $\left\lbrace \rho_{\tau}^k \right\rbrace_{k \in \mathbb{N}}$ para fabricar funciones de $\left[ 0, \infty \right)$ a $L^1 \left( \mathbb{R}^2 \right)$. Para cada entero positivo $k$, sea $\nabla \varphi^k$ el optimal transportation plan con $\nabla \varphi^k \sharp \rho_{\tau}^{k+1} = \rho_{\tau}^k$ (puede consultarse el Apéndice \ref{ap:Transporte_optimo} para encontrar las definiciones de estos objetos). Entonces, para $k \tau \leq t \leq (k+1) \tau$ definimos
\begin{equation*}
\rho_{\tau} (t) = \left( \dfrac{t - k\tau}{\tau} \mathrm{id} + \dfrac{(k+1) \tau - t}{\tau} \nabla \varphi^k \right) \sharp \rho_{\tau}^{k+1} \, .
\end{equation*}
Usando esta interpolación y el Teorema de Benamou-Brenier (consúltese el Apéndice \ref{ap:Transporte_optimo}), se observa que satisface $\rho_{\tau} (k \tau ) = \rho_{\tau}^k$, $\rho_{\tau} ((k+1) \tau ) = \rho_{\tau}^{k+1}$ y $\mathcal{W}_2 ( \rho_{\tau}^k, \rho_{\tau} (t)) = (t-k \tau) \mathcal{W}_2 (\rho_{\tau}^k , \rho_{\tau}^{k+1} )$.

Finalmente, para poder cumplir los objetivos que nos habíamos marcado al comienzo del Capítulo, necesitamos ser capaces de aplicar toda la teoría que hemos desarrollado durante las Secciones \ref{sec:Sub_energia_libre} y \ref{sec:Sub_JKO} al caso continuo, es decir, comprobar que la regla \eqref{ec:JKO} converge cuando $\tau \rightarrow 0$ \cite{BLA_CAL_CAR},
\begin{teor}[Convergencia del método \eqref{ec:JKO} cuando $\tau \rightarrow 0$] \label{th:Masa_subcritica_JKO}
Si $M< 8 \pi$, entonces la familia $\left( \rho_{\tau} \right)_{\tau > 0}$ admite una sub-sucesión que converge débilmente en $L^1 \left( \mathbb{R}^2 \right)$ a una solución débil del sistema \eqref{ec:KS_clasico}: para todo $\left( t_1, t_2 \right) \subseteq \left[ 0, + \infty \right)$, y para todo $\zeta$ suave
\begin{equation}\label{ec:Teorema_Sub}
\begin{array}{ll}
\dfrac{d}{dt} {\displaystyle \int}_{\mathbb{R}^2} \zeta (x) \rho (x,t) \, dx & = \int_{\mathbb{R}^2} \Delta \zeta (x) \rho (x, s) \, dx \, ds \\
 & - \dfrac{1}{4 \pi} {\displaystyle \int} {\displaystyle \int}_{\mathbb{R}^2 \times \mathbb{R}^2} \rho (x, s) \rho (y, s) \dfrac{(x-y) \cdot (\nabla \zeta (x) - \nabla \zeta (y))}{| x- y |^2} \, dy \, dx .
\end{array}
\end{equation}
\end{teor}

\section{Demostración del Teorema \ref{th:Masa_subcritica_JKO}}\label{sec:Sub_dem}
Es estudiado por Blanchet, Calvez y Carrillo en \cite{BLA_CAL_CAR}, siguiendo las líneas generales de la prueba de convergencia del método para el flujo gradiente euclídeo. Nosotros solo presentamos las ideas principales de la misma, siguiendo las notas de Blanchet \cite{BLA}.

\textit{(i) \textbf{Existencia de minimizadores}:} Hay que enfatizar que el funcional $E_{KS}$ no es convexo, ya que tiene un potencial de interacción no-convexo. Por ello, incluso la existencia de un minimizador no está claro, y por ese mismo motivo tampoco se puede decir nada sobre la unicidad. No obstante, podemos construir una sucesión de minimizadores de \eqref{ec:JKO} cuando $M < 8 \pi$ si recurrimos a la estimación \eqref{ec:Subcritico_Estimacion}.

\textit{(ii) \textbf{La ecuación discreta de Euler-Lagrange para el funcional de \eqref{ec:JKO}}:} La perturbación del minimizador se ha hecho utilizando el transporte óptimo: Sea $\zeta$ un campo vectorial suave con soporte compacto, introducimos $\psi_{ \varepsilon} := |x|^2 / 2 + \varepsilon \zeta$. Definimos $\bar{\rho_{\varepsilon}}$, la perturbación push-forward de $\rho_{\tau}^{n+1}$ por $\nabla \psi_{\varepsilon}$:
\begin{equation*}
\bar{\rho_{\varepsilon}} = \nabla \psi_{\varepsilon} \sharp \rho_{\tau}^{n+1} .
\end{equation*}
Tras algunos cálculos estándar, obtenemos
\begin{equation*}
\begin{array}{ll}
{\displaystyle \int}_{\mathbb{R}^2} \nabla \zeta (x) \dfrac{x - \nabla \varphi^n (x)}{\tau} \rho_{\tau}^{n+1} (x) \, dx \\
\, \, \, \, \, \, \, \, \, \, \, \, \, \, \, \, \, \, \, \, \, \, \, \, \, \, \, \, \, \, = {\displaystyle \int}_{\mathbb{R}^2} \left[ \Delta \zeta (x) - \dfrac{1}{4 \pi} {\displaystyle \int}_{\mathbb{R}^2} \dfrac{\left[ \nabla \zeta (x) - \nabla \zeta (y) \right] \cdot (x-y)}{|x-y|^2} \rho_{\tau}^{n+1} (y) \, dy \right] \rho_{\tau}^{n+1} (x) \, dx ,
\end{array}
\end{equation*}
que es la forma débil de la ecuación de Euler-Lagrange:
\begin{equation}
\dfrac{id - \nabla \varphi^n}{\tau} \rho_{\tau}^{n+1} = - \nabla \rho_{\tau}^{n+1} + \rho_{\tau}^{n+1} \nabla c_{\tau}^{n+1}  .
\end{equation}
Seguidamente, usamos la expansión de Taylor para comprobar que hemos dado bien el paso del discreto al continuo
\begin{equation*}
\zeta (x) - \zeta \left[ \nabla \varphi^n (x) \right] = \left[ x - \nabla \varphi^n (x) \right] \cdot \nabla \zeta (x) + O \left[ \left| x- \nabla \varphi^n (x) \right|^2 \right] ,
\end{equation*}
obteniendo, para todo $t_2 > t_1 \geq 0$,
\begin{equation}\label{ec:subcritico_formulacion_debil}
\begin{array}{ll}
{\displaystyle \int}_{\mathbb{R}^2} \zeta (x) \left[ \rho_{\tau} (x, t_2) - \rho_{\tau} (x, t_1) \right] \, dx  = {\displaystyle \int}_{t_1}^{t_2} {\displaystyle \int}_{\mathbb{R}^2} \Delta \zeta (x) \rho_{\tau} (x, s) \, dx \, ds + O (\tau^{1/2} ) \\
 \, \, \, \, \, \, \, \, \, \, \, \, \, \, \, \, \, \, \, \, \, \, \, \, \, \, \, \, \, \, \, \, \, \, \, \, \, \, \, \, \, \, \, \, \, \,   - \dfrac{1}{4 \pi} {\displaystyle \int}_{t_1}^{t_2} {\displaystyle \int} {\displaystyle \int}_{\mathbb{R}^2 \times \mathbb{R}^2} \rho_{\tau} (x, s) \rho_{\tau} (y, s) \dfrac{(x-y) \cdot \left( \nabla \zeta (x) - \nabla \zeta (y) \right)}{|x-y|^2} \, dy \, dx ,
\end{array}
\end{equation}
y recuperamos \eqref{ec:Teorema_Sub}.

\textit{(iii)} \textbf{\textit{Una estimación} a priori\textit{:}} Para pasar al límite la regla \eqref{ec:JKO} nos da algunas acotaciones \textit{a priori}. Tenemos:
\begin{equation}\label{ec:subcritica_estimacion_a_priori}
E_{KS} \left[ \rho_{\tau}^{n+1} \right] + \dfrac{1}{2 \tau} \mathcal{W}_2^2 \left( \rho_{\tau}^n, \rho_{\tau}^{n+1} \right) \leq E_{KS} \left[ \rho_{\tau}^n \right] .
\end{equation}
Como consecuencia obtenemos una \textit{estimación de la energía}
\begin{equation}\label{ec:subcritica_energy_estimate}
\sup_{n \in \mathbb{N}} E_{KS} \left[ \rho_{\tau}^n \right] \leq E_{KS} \left[ \rho_{\tau}^0 \right]
\end{equation}
y una \textit{estimación total cuadrática}
\begin{equation}\label{ec:subcritico_Holder_continuo}
\dfrac{1}{2 \tau} \sum_{n \in \mathbb{N}} \mathcal{W}_2^2 \left( \rho_{\tau}^n , \rho_{\tau}^{n+1} \right) \leq E_{KS} \left[ \rho_{\tau}^0 \right] - \inf_{n \in \mathbb{N} } E_{KS} \left[ \rho_{\tau}^n \right] .
\end{equation}

\textit{(iv) \textbf{Pasando al límite}:} La estimación de la energía \eqref{ec:subcritica_energy_estimate} y \eqref{ec:Subcritico_Cota_Energia_Libre} nos da una cota para $\int \rho \log \rho$, al menos mientras $M < 8 \pi$. La acotación de $\rho_{\tau} \log \rho_{\tau}$ demuestra que la solución no tiene un blow-up. De hecho, usando
\begin{equation*}
\int_{\left\lbrace \rho > K \right\rbrace } \rho \, dx \leq \dfrac{1}{\log K} \int_{\left\lbrace \rho > K \right\rbrace } \rho | \log \rho | \, dx \leq \dfrac{C}{\log K},
\end{equation*}
obtenemos que $\left( \rho_{\tau} \right)_{\tau}$ converge a una determinada $\rho$ en $weak-L^1 \left( \mathbb{R}^2 \right)$. Podemos usar la $1/2$-Hölder continuidad \eqref{ec:subcritico_Holder_continuo} y el Teorema de Ascoli para obtener convergencia en $C^0 \left( \left[ 0, T \right] ; \mathcal{P} \left( \mathbb{R}^2 \right) \right)$. 

De esta manera, podemos pasar al límite $\tau \rightarrow 0$ en \eqref{ec:subcritico_formulacion_debil} y demostrar que $\rho$ es una solución débil de \eqref{ec:flujo_Wasserstein}. Obsérvese que el último término de \eqref{ec:subcritico_formulacion_debil} converge porque la convergencia de $\left( \rho_{\tau} \right)_{\tau} $ en $weak-L^1 \left( \mathbb{R}^2 \right)$ asegura convergencia de $\left( \rho_{\tau} \bigotimes \rho_{\tau} \right)_{\tau}$ en $weak-L^1 \left( \mathbb{R}^2 \right)$. De todas formas, las soluciones construidas durante este Capítulo son débiles.

\chapter{Masa crítica}\label{ch:Critica}

En este Capítulo estudiaremos el caso de la masa crítica, el más delicado de los tres y que se alcanza cuando $M$ toma el valor $8 \pi$. Bajo esta hipótesis, por el trabajo de Blanchet, Carrillo y Masmoudi \cite{BLA_CAR_MAS}, sabemos que se preserva el segundo momento de la solución con el tiempo y que las soluciones están definidas globalmente en tiempo. Si el segundo momento es infinito, existe una gran variedad de comportamientos observados por ejemplo en \cite{LOP_NAG_YAM, LOP_NAG_YAM_2, NAG_YAM}, exhibiendo, entre otros, complejos comportamientos oscilatorios. En cambio, si el segundo momento inicial es finito, existe blow-up en tiempo infinito para las soluciones, \cite{BLA_CAR_MAS}. Nosotros nos centraremos en el segundo caso, siguiendo los resultados de Dávila, del Pino, Dolbeault, Musso y Wei \cite{DAV_PIN_DOL_MUS_WEI}. Nuestro objetivo principal será comprender cuál es la forma que toman estos blow-up. El siguiente Teorema es el principal resultado del Capítulo, enunciado por primera vez en \cite{PIN} y demostrado en \cite{DAV_PIN_DOL_MUS_WEI}. La idea principal es que los blow-up toman la forma de una bubble.
\begin{teor}\label{th:masa_critica}
Existe una función con simetría radial no-negativa $\rho_0^{\ast} \left( x \right)$ con masa crítica $\int_{\mathbb{R}^2} \rho_0^{\ast} ( x ) \, dx = 8 \pi$ y segundo momento finito $\int_{\mathbb{R}^2} |x|^2 \rho_0^{\ast} (x) \, dx < + \infty$ y tal que para cada dato inicial $\rho_0$ suficientemente cerca (en el sentido correcto) a $\rho_0^{\ast}$ con $\int_{\mathbb{R}^2} \rho_0 \, dx = 8 \pi$, tenemos que la solución $\rho (x,t)$ asociada al sistema \eqref{ec:KS_clasico} tiene la forma
\begin{equation*}
\rho (x, t) = \dfrac{1}{\lambda (t)^2} U_0 \left( \dfrac{x - \xi (t)}{\lambda (t)} \right) (1 + o (1) ), \, \, \, \, \, \, \, \, \, \, \, \, \, \, \, \, \, \, U_0 (y) = \dfrac{8}{( 1 + | y |^2 )^2}
\end{equation*}
uniformemente en conjuntos acotados de $\mathbb{R}^2$, y 
\begin{equation}\label{ec:Critico_lambda}
\lambda (t) = \frac{k}{\sqrt{\log t}} (1 +o(1)), \, \, \, \xi (t) \rightarrow q, \, \, \, \mathrm{cuando} \, \, \, t \rightarrow + \infty
\end{equation}
para algún número $k > 0$ y alguna $q \in \mathbb{R}^2$.
\end{teor} 
En este contexto, suficientemente cerca significa $\rho_0 (x) := \rho_0^{\ast} (x) + \varphi (x)$ con 
\begin{equation*}
\| \varphi \|_{\ast \sigma} := \| (1 + | \cdot |^{4+\sigma}) \varphi \|_{L^{\infty} (\mathbb{R}^2)} + \| (1 + | \cdot |^{5+ \sigma} \nabla \varphi (x) \|_{L^{\infty} (\mathbb{R}^2)} < + \infty
\end{equation*}
para algún $\sigma >0$.

No vamos a dar una demostración completa del Teorema. Nos conformaremos con una derivación formal de $\lambda$. Durante esta derivación formal sí que demostraremos que el blow-up tiene lugar en el centro de masas de la densidad celular, un fenómeno que fue conjeturado por primera vez por Childress y Percuss \cite{CHI_PER}. Para el lector interesado, se puede acudir a unas notas de Herrero y Velázquez \cite{HER_VEL}.

Comenzamos señalando que sabemos que la función $U_0$ del Teorema pertenece a una familia de soluciones estacionarias de \eqref{ec:KS_clasico} para el caso crítico $M= 8 \pi$.
\begin{equation}\label{ec:Burbujas}
U_{\lambda , \xi} \left( x \right) = \dfrac{1}{\lambda^2} U_0 \left( \dfrac{x - \xi}{\lambda} \right) , \, \, \, U_0 \left( y \right) = \dfrac{8}{( 1 + | y |^2 )^2}, \, \, \, \lambda >0, \, \xi \in \mathbb{R}^2 .
\end{equation}
De esto deducimos que la familia de funciones $\left( U_{\lambda , \xi} \right)$ con forma de bubbles, son de estados estacionarios de masa finita (crítica) y segundo momento infinito, y forman un conjunto de soluciones de \eqref{ec:KS_clasico}. Si una solución de \eqref{ec:KS_clasico} es atraída por la familia $ U_{ \lambda , \xi} $, su masa no puede ser inferior a $8 \pi$ y tener segundo momento finito. De ser así (consúltese Capítulo \ref{ch:Subcritica}), el segundo momento permanecería finito a lo largo de las soluciones $\forall \, t \geq 0$ por lo que no se produciría blow-up, contradictorio con aproximarse a la familia $U_{\lambda , \xi}$, donde este sí tiene lugar.  

Recordemos que para nuestro caso de interés, masa crítica $M=8 \pi$ y segundo momento finito, el blow-up ocurre en tiempo infinito. En él, la solución toma la forma de burbuja descrita en \eqref{ec:Burbujas} con $\lambda = \lambda \left( t \right) \rightarrow 0$. Vamos a calcular formalmente quién es $\lambda$.

Tomamos como primera solución una $\rho \left( x, t \right)$ genérica globalmente definida en tiempo con masa crítica y segundo momento finito consistente con un comportamiento de la forma 
\begin{equation}\label{ec:Burbuja_modelo_u}
\rho \left( x, t \right) = \dfrac{1}{\lambda \left( t \right)^2} U_0 \left( \dfrac{x - \xi \left( t \right)}{ \lambda \left( t \right)} \right) \left( 1 + o \left( 1 \right) \right) \, \, \mathrm{cuando} \, \, t \rightarrow \infty
\end{equation}
para algunas funciones determinadas $0 < \lambda \left( t \right) \rightarrow 0$ y $\xi \left( t \right) \rightarrow q \in \mathbb{R}^2$. 

Debido a que el segundo momento de $U_0$ es infinito, no esperamos que la aproximación \eqref{ec:Burbuja_modelo_u} sea uniforme en $\mathbb{R}^2$. No obstante, suficientemente lejos, sí que observamos un decaimiento rápido en $x$ para el que somos capaces de encontrar una expresión de la aproximación asintótica para el parámetro de reescalamiento $\lambda \left( t \right)$ que encaje con este comportamiento.

Para ello introducimos la función $C_0 := \left( - \Delta \right)^{-1} U_0$. Realizando algunos cálculos vemos que $C_0 \left( y \right) = \log \frac{8}{( 1 + | y |^2 )^2}$. También nos damos cuenta de que resuelve la ecuación de Liouville (véase Apéndice \ref{ap:Liouville}) 
\begin{equation}\label{ec:Ecuacion_Liouville}
- \Delta C_0 = e^{C_0} = U_0 \, \, \mathrm{en} \, \, \mathbb{R}^2.
\end{equation}
Entonces $\nabla C_0 \left( y \right) \approx - \frac{4 y}{ | y|^2}$ para todo $y$ grande. Lejos de $x = \xi$ se satisface, $- \nabla \cdot ( \rho \nabla ( - \Delta )^{-1} \rho ) \approx 4 \nabla \rho \cdot \frac{x - \xi}{| x - \xi |^2}$. Usamos coordenadas polares para escribir nuestras expresiones, $\rho ( r, \theta, t ) = \rho ( x, t )$, $x = \xi ( t ) + r e^{i \theta}$. Esto y asumir $\dot{\xi} ( t ) \rightarrow 0$ suficientemente rápido, implica que la ecuación \eqref{ec:KS_clasico} se lee aproximadamente de la forma $\partial_t \rho = \partial_r^2 \rho + \frac{5}{r} \partial_r \rho$.

Por tanto, se puede usar la ecuación del calor homogénea en $\mathbb{R}^6$ para funciones con simetría radial para entender esta función. Para nuestra ecuación la región auto-similar es $r \simeq \sqrt{t}$. Lejos de la región auto-similar, $r \gg \sqrt{t}$ la función que estamos estudiando presenta decaimiento rápido en $+ \infty$, como consecuencia de la finitud del segundo momento. 

Para los valores pequeños, $r \ll \sqrt{t}$ no conocemos su comportamiento y debemos estudiarlo de manera más exhaustiva. El método que usamos para estudiar esta región consiste en truncar la burbuja \eqref{ec:Burbujas} más allá de la zona auto-similar. Esto se puede conseguir usando alguna función cut-off. Introducimos el parámetro $\alpha$ y definimos
\begin{equation}\label{ec:Burbujas_u1}
\rho_1 \left( x, t \right) = \dfrac{\alpha \left( t \right)}{\lambda^2} U_0 \left( \dfrac{| x - \xi |}{\lambda} \right) \chi \left( x, t \right),
\end{equation}
denotando $\chi \left( x, t \right) = \chi_0 \left( \frac{| x - \xi |}{\sqrt{t}} \right)$ con $\chi_0 \left( s \right)$ una función cut-off suave que toma el valor $1$ para $s \leq 1$ y $0$ para $s \geq 2$. El otro parámetro que hemos introducido, $\alpha ( t )$, sirve para poder imponer la condición $\int_{\mathbb{R}^2} \rho_1 \left( x, t \right) dx = 8 \pi$, $\forall \, t \geq 0$. Realizando algunos cálculos directos deducimos que para nuestro caso  
\begin{equation*}
\alpha ( t ) = 1 + a \frac{\lambda^2}{t} \left( 1 + o ( 1 ) \right), \, \, \mathrm{con} \, \, a = 2 \int_0^{\infty} \frac{1- \chi_0 ( s )}{s^3} ds.
\end{equation*}
Queremos que, aproximadamente, el valor de reescalamiento $\lambda (t)$ sea de tal forma que este tenga consistencia con el empleo de una solución que satisfaga $\rho ( x, t ) \approx \rho_1 ( x, t )$. Como $\rho ( x, t )$ es solución de \eqref{ec:KS_clasico} cumple
\begin{equation}\label{ec:Burbujas_Conservaciones}
\dfrac{d}{dt} \int_{\mathbb{R}^2} \rho \left( x, t \right) \left| x \right|^2 dx = 4M - \dfrac{M^2}{2 \pi}, \, \, \, M = \int_{\mathbb{R}^2} \rho \left( x, t \right) \, dx .
\end{equation}
Además, los primeros momentos de $\rho$ (centro de masa de $\rho$ en espacio) se conservan, 
\begin{equation}\label{ec:Critico_centro_de_masas}
\frac{d}{dt} \int_{\mathbb{R}^2} \rho ( x, t ) x_i \, dx =0.
\end{equation}
La solución $\rho \left( x, t \right)$ toma valores cercanos a los de $\rho_1 \left( x, t \right)$. Como $ \int_{\mathbb{R}^2} \rho_1 ( x, t ) \, dx = 8 \pi$ es posible deducir que la identidad $\frac{d}{dt} \int_{\mathbb{R}^2} \rho_1 \left( x, t \right) \left| x \right|^2 \, dx =0$ es cierta, al menos de forma aproximada. Usamos una expresión más clara para nuestros intereses,
\begin{equation*}
a I \left( t \right) := \int_{\mathbb{R}^2} \dfrac{\alpha_0}{\lambda^2} U_0 \left( \dfrac{x - \xi}{\lambda} \right) \chi_0 \left( \dfrac{ | x - \xi |}{\sqrt{t}} \right) | x |^2 dx \, = \, \mathrm{constante}, 
\end{equation*}
en la que es fácil comprobar que para alguna constante $\kappa$ 
\begin{equation*}
I \left( t \right) = 16 \pi \lambda^2 \int_0^{\sqrt{t} / \lambda} \frac{u^3 d u}{( 1 + u^2 )^2} + \kappa + o \left( 1 \right) = 16 \pi \lambda^2 \log \frac{\sqrt{t}}{\lambda} + \kappa + o \left( 1 \right). 
\end{equation*}
De todo esto se sigue que $\lambda ( t )$ satisface aproximadamente
$\lambda^2 \log t \, = \, k^2 \, = \, \mathrm{constante}$, es decir, $\lambda (t)$ cumple \eqref{ec:Critico_lambda}, y hemos acabado la demostración. 

Por último, y gracias a \eqref{ec:Critico_centro_de_masas}, nos damos cuenta de que el centro de masas se preserva a lo largo de una solución verdadera, $\frac{d}{dt} \int_{\mathbb{R}^2} x \rho ( x, t ) dx = 0$, que demuestra la conjetura de Childress y Percuss que habíamos señalado al comienzo del capítulo. 

Además, como el centro de masas de $\rho_1 (x, t)$ es exactamente $\xi (t)$ estimamos su valor aproximado, $\xi (t) \, = \, \mathrm{constante} \, = \, q $, obteniendo el último detalle para comprender por completo la familia \eqref{ec:Burbujas}.



\chapter{Masa Supercrítica}\label{ch:Supercritica}
Tal y como señalamos en la Introducción, para masa $M > 8 \pi$, si empezamos con segundo momento finito, este se vuelve negativo a partir de un cierto $T$, esto es una contradicción que surge como consecuencia de suponer que la solución esta definida para todo tiempo $t>0$. Por lo tanto, debe existir blow-up en tiempo finito y la agregación domina sobre la difusión. Un blow-up puede ser de tipo I o de tipo II (consúltese el Apéndice \ref{ap:auto-similar}). Decimos que una solución $\rho (t)$ de \eqref{ec:KS_clasico} exhibe blow-up tipo I en tiempo $t=T$ si existe una constante $C>0$ tal que
\begin{equation*}
{\lim \sup}_{t\rightarrow T} (T-t) \| \rho (t) \|_{L^{\infty} (\mathbb{R}^2)} \leq C. 
\end{equation*}
Las singularidades que estudiamos en este capítulo son todas tipo II. 

En el trabajo pionero \cite{HER_VEL}, Herrero y Velázquez usan \textit{matching asymptotics} para conseguir el primer ejemplo de una solución radial en la que se produce blow-up en tiempo finito, con un estudio posterior más detallado de la estabilidad hecho por Velázquez en \cite{VEL, VEL_2, VEL_3}. Más adelante, en \cite{RAP_SCH}, Raphaël y Schweyer demuestran la existencia de la formación de singularidades en el caso radial con masa cercana a $8 \pi$, obteniendo una descripción completa del proceso de agregación asociado. Para hacer esto, se recurre al estudio del espectro del operador linealizado y a la construcción de familias aproximadas de perfiles auto-similares. Recientemente, refinando el estudio del espectro, Collot, Ghoul, Masmoudi y Nguyen en \cite{COL_GHO_MAS_NGU}, han conseguido eliminar la hipótesis de masa inicial cercana a $8 \pi$ necesaria en el trabajo anterior. Proporcionando una nueva aproximación a las singularidades tipo II para problemas parabólicos críticos, con técnicas que pueden adaptarse a más tipos de ecuaciones con blow-up.

Aquí solo construiremos el ejemplo de Herrero y Velázquez \cite{HER_VEL} porque contiene las ideas fundamentales para el resto de artículos posteriores. Queremos mostrar una clase de soluciones de \eqref{ec:KS_clasico} para las cuales ocurre un colapso quimiotáctico y describirlo con precisión. Para construirlo, vamos a añadir una condición a \eqref{ec:KS_clasico}, se define solo para $x \in \Omega$, con $\Omega$ un conjunto abierto, regular y acotado. Además, cambiamos la condición de frontera que enunciamos en \eqref{ec:KS_clasico} por otra de no-flujo (la población de bacterias se mantiene constante)
\begin{equation}\label{ec:Neumann}
\frac{\partial \rho}{\partial n} = \frac{\partial c}{\partial n} = 0 \, \, \, \, \, \mathrm{para} \, \, \, \, \, x \in \partial \Omega, \, \, \, \, \,  t > 0.
\end{equation} 
En particular, se obtiene el siguiente:

\begin{teor}[Herrero-Velázquez, \cite{HER_VEL}]\label{th:Herrero_Velazquez_1}
Consideremos un sistema que satisface \eqref{ec:KS_clasico} (definida solo para $x \in \Omega$), \eqref{ec:Neumann}, para tiempos suficientemente pequeños $T>0$ en una bola $\Omega \equiv B_R = \left\lbrace x : |x| < R \right\rbrace $, con $R$ cualquier número positivo dado. Entonces, para cualquier $T > 0$ fijo existe una solución radial $\rho (r, t)$ que experimenta un blow-up para $r=0$ y $t=T$ y tal que:
\begin{equation}\label{ec:Super_T1_rho}
\rho (r, t) \rightarrow 8 \pi \delta (x) + f (r) \, \, \, \, \, \mathrm{cuando} \, \, \, \, \, t \uparrow T
\end{equation}
en el sentido de distribución, donde
\begin{equation}\label{ec:Super_T1_f}
f(r) = \frac{C}{r^2} e^{- 2 | \log r |^{1/2}} (2 | \log r | )^{\frac{1}{2 \sqrt{2} | \log r |^{1/2}} - \frac{1}{2}} (1 + o (1) )
\end{equation}
cuando $r \rightarrow 0$, y $C$ una constante positiva. Para tiempo $t = T$, el perfil cerca de $r=0$ viene dado por
\begin{equation}
\rho (r, T) = 8 \pi \delta (x) + f(r)
\end{equation}
con $f$ como en \eqref{ec:Super_T1_f}. Por último, si definimos
\begin{equation}
S(t) = (T-t) ( \sup_{\Omega} \rho (r, t) ) \equiv (T-t) \rho (0, t)
\end{equation}
uno tiene que $\lim_{t \uparrow T} S(t) = \infty$. En concreto, 
\begin{equation}
S(t) \sim C_1 (T- t)^{-1/2} | \log (T-t) |^{1-\frac{1}{\sqrt{| \log (T -t) |}}} \, \, \, \, \, \mathrm{cuando} \, \, \, \, \, t \uparrow T
\end{equation}
para alguna constante $C_1 > 0$.
\end{teor}
Hay que destacar que las soluciones a las que nos referimos en el Teorema \ref{th:Herrero_Velazquez_1} están definidas en el intervalo temporal $(T - \delta , T)$, con $\delta >0 $ suficientemente pequeño. También es interesante resaltar que en este ejemplo podemos estudiar como se forma el blow-up $\delta (x)$ cerca de $T$. Esto queda recogido en:
\begin{teor}[Herrero-Velázquez \cite{HER_VEL}]\label{th:Herrero_Velazquez_2}
Sea $\rho (r, t)$ la solución descrita en el Teorema \ref{th:Herrero_Velazquez_1}. Entonces, uno tiene que
\begin{equation}\label{ec:HV2_forma_rho}
\rho (r, t) = \left( \frac{1}{R(t)^2} U_0 \left( \frac{r}{R(t)} \right) \right) (1 + o(1)) + o \left( \frac{e^{- \sqrt{2} | \log (T-t)|^{1/2}}}{r^2} \cdot \mathbbm{1}_{\left\lbrace r \geq R(t) \right\rbrace } \right)
\end{equation}
cuando $t\uparrow T$, uniformemente en conjuntos $r \leq C (T-t)^{1/2}$, con $U_0(r) = \frac{8}{(1+r^2)^2}$ el estado estacionario, y, también cuando $t \uparrow T$
\begin{equation}\label{ec:HVT2_perfil_R}
R(t) = C (T-t)^{1/2} e^{- \sqrt{2} / 2 | \log (T-t) |^{1/2}} ( | \log (T-t) | )^{1/4 ( | \log (T-t) |)^{-1/2} - 1/4} (1 + o(1)).
\end{equation}
\end{teor}

La demostración de estos dos Teoremas se basa en el método de las \textit{matched asymptotic expansions} que desarrollamos en la Sección \ref{sec:matched_asymptotic}. No damos la prueba completa, solo el cálculo formal de los asymptotics. Aunque, antes de comenzar con estas ideas, en la Sección \ref{sec:transformacion_integral}, vamos a transformar el problema en uno más adecuado a nuestros objetivos. Como ya hemos dicho al comienzo de este Capítulo seguiremos el paper clásico \cite{HER_VEL}.

\section{Una transformación integral}\label{sec:transformacion_integral}
Sea $\rho (r, t)$ una solución radial de las ecuaciones de Keller-Segel con condiciones de frontera Neumann, definida en una bola $B_R \subseteq \mathbb{R}^2$ para tiempos $0 < t < T < \infty$. Definimos $m (r, t)$ como
\begin{equation*}
m(r, t) = \int_{\left\lbrace  |x| \leq r \right\rbrace } ( \rho (x, t) - 1 ) \, dx \, \, \, \, \mathrm{si} \, \, \, \, 0 < r < R, \, \, \, \, \, \, m(r, t) = 0 \, \, \, \, \mathrm{si} \, \, \, \, r > R.
\end{equation*}
Obsérvese que la condición de frontera $\frac{\partial c}{\partial n} = 0$ en $r=R$ implica que $m(R,t) = 0$. 

Tras definir $m$, damos el siguiente resultado, muy importante para la aproximación a la demostración.
\begin{lema}
Sea $R>0$ y $T>0$, y sea $\rho (r, t)$ una solución radial para $r < R$ y $t < T$. Consideremos la función $w (r,t)$ para el centro de masas
\begin{equation}\label{ec:Super_Lema_w}
w (r,t) = \int_0^r \xi m ( \xi , t ) \, d \xi
\end{equation}
Entonces, si introducimos variables auto-similares dadas por
\begin{equation}\label{ec:Super_Lema_auto-similares}
W (y , \tau ) = (T- t)^{-1} w (r, t), \, \, \, \, \mathrm{donde} \, \, \, \, y = r (T-t)^{-1/2}, \, \, \, \, \tau = - \log (T-t),
\end{equation}
resulta que $W (y, \tau )$ satisface la ecuación:
\begin{equation}\label{ec:Super_Lema_W}
W_{\tau} = W_{yy} + \left( \frac{1}{y} - \frac{y}{2} \right) W_y + W + \left( \frac{1}{4 \pi} \left( \frac{W_y}{y} \right)^2 - 4 \left( \frac{W_y}{y} \right) \right) + e^{- \tau} W
\end{equation}
cuando $0 < y < \infty$ y $ \tau > 0$.
\end{lema}
\begin{proof}
Esta prueba consiste en un cálculo directo. Una comprobación rápida revela que $m(r,t)$ satisface 
\begin{equation}\label{ec:Super_Lema_dem}
m_t = m_{rr} + \left( \frac{m}{2 \pi} -1 \right) \frac{m_r}{r} + m.
\end{equation}
Se sigue de \eqref{ec:Super_Lema_w} y \eqref{ec:Super_Lema_dem} que $w (r, t)$ resuelve
\begin{equation*}
w_t = w_{rr} + \frac{w_r}{r} + \left( \frac{1}{4 \pi} \left( \frac{w_r}{r} \right)^2 - 4 \left( \frac{w_r}{r} \right) \right) + w
\end{equation*} 
para $r > 0$ y $0 < t < T$, donde \eqref{ec:Super_Lema_W} se sigue de un cambio de variables de \eqref{ec:Super_Lema_auto-similares}.
\end{proof}
En lo siguiente, también recurriremos a la función auxiliar $\phi (y, \tau )$ dada por la relación 
\begin{equation}\label{ec:Super_def_phi}
W(y, \tau ) = \int_0^y r \phi (r, \tau ) \, dr
\end{equation}
que satisface la ecuación
\begin{equation}\label{ec:Super_ec_en_phi}
\phi_{\tau} = \phi_{yy} - \frac{1}{2} y \phi_y + \left( \frac{1}{2 \pi} \phi - 1 \right) \frac{\phi_y}{y} + e^{- \tau} \phi.
\end{equation}
Obsérvese que las funciones $m$ y $\phi$ satisfacen $m(r,t) = \phi (y, \tau )$, con $y$, $\tau$ las variables autosimilares dadas en \eqref{ec:Super_Lema_auto-similares}.

\section{El enfoque de matched asymptotic expansions}\label{sec:matched_asymptotic}
Para comenzar, suponemos que $W( y, \tau )$ se aproximará al estado estacionario de \eqref{ec:Super_Lema_W} cuando $\tau \rightarrow \infty$. Un argumento de contar órdenes sugiere entonces que
\begin{equation}\label{ec:Super_estado_estacionario}
W (y, \tau ) \rightarrow 4 \pi y^2 \, \, \, \, \mathrm{cuando} \, \, \, \, \tau \rightarrow \infty,
\end{equation}
que se consigue inmediatamente al suponer $W_{\tau} \to 0$, $e^{-\tau} W \to 0$ cuando $\tau \to \infty$. Uno espera que \eqref{ec:Super_estado_estacionario} se satisfaga uniformemente en conjuntos compactos de $y$. De \eqref{ec:Super_def_phi} y \eqref{ec:Super_estado_estacionario} se deduce que $\phi (y, \tau ) \rightarrow 8 \pi$ cuando $\tau \rightarrow \infty$. Por otra parte, $\phi (0, \tau ) = 0$ para cualquier $\tau > 0$ (cf. \eqref{ec:Super_estado_estacionario}), lo que nos guía a pensar en el inicio de una \textit{boundary layer} cerca del origen donde cambios abruptos deberían aparecer cuando $\tau \rightarrow \infty$. Para estudiar la boundary layer denotamos por $\varepsilon ( \tau )$ su tamaño. Entonces, introducimos una nueva escala espacial para la región interior por medio del ajuste
\begin{equation}\label{ec:Super_xi}
\xi = \tfrac{y}{\varepsilon (\tau )} .
\end{equation}
Si recordamos \eqref{ec:Super_ec_en_phi}, ya hemos visto que $\phi ( \xi , \tau )$ satisface
\begin{equation}\label{ec:Super_ec_phi_xi}
\varepsilon^2 \phi_{\tau} - \varepsilon \dot{\varepsilon} \xi \phi_{\xi } = \phi_{ \xi \xi } + \left( \frac{\phi}{2 \pi} - 1 \right) \frac{\phi_{\xi}}{\xi} - \varepsilon^2 \frac{\xi \phi_{\xi}}{2} + \varepsilon^2 e^{ - \tau} \phi,
\end{equation}
y como estamos asumiendo $\lim_{\tau \rightarrow \infty} \varepsilon ( \tau ) = 0$, la ecuación \eqref{ec:Super_ec_phi_xi} sugiere que, para $\tau \gg 1$, $\phi (y, \tau )$ debería  aproximar una solución de la ecuación estacionaria
$\xi \phi_{\xi \xi}  + \left( \frac{\phi}{2 \pi} - 1 \right) \phi_{\xi} = 0$, con $\phi (0) = 0$, y cuyas soluciones vienen dadas por
\begin{equation}\label{ec:Super_estados_estacionarios_phi_a}
\phi_a ( \xi ) = 8 \pi \cdot \frac{\xi^2}{\xi^2 + a^2},
\end{equation}
donde $a$ es un número real arbitrario. En vista de \eqref{ec:Super_xi}, podemos normalizar la elección de $\epsilon (\tau )$ para que tengamos $a = 1$ en \eqref{ec:Super_estados_estacionarios_phi_a}. Es decir, esperamos el siguiente comportamiento asintótico para $\phi$
\begin{equation}\label{ec:Super_estados_estacionarios_phi_1}
\phi (y, \tau ) \sim 8 \pi \cdot \frac{y^2}{y^2 + \varepsilon ( \tau )^2} \, \, \, \, \mathrm{cuando} \, \, \, \, \tau \rightarrow \infty .
\end{equation}
Tras esto, lo siguiente es determinar el valor de $\epsilon ( \tau )$. Para lograrlo hacemos uso de la función $W (y, \tau )$, teniendo en cuenta \eqref{ec:Super_estado_estacionario} podemos escribir $W (y, \tau ) = \psi (y, \tau ) + 4 \pi y^2$ para estudiar cómo se acerca $W$ a su estado estacionario. Por \eqref{ec:Super_def_phi} y \eqref{ec:Super_Lema_W}, $\psi$ resuelve 
\begin{equation}\label{ec:Super_ec_psi}
\psi_{\tau } = \psi_{yy} + \left( \frac{1}{y} - \frac{y}{2} \right) \psi_y + \psi + \left( \frac{1}{4 \pi} \phi^2 - 4 \phi + 16 \pi \right) + e^{- \tau} \phi.
\end{equation}
Con esto, si \eqref{ec:Super_estados_estacionarios_phi_1} se satisface, entonces se puede asegurar que se tiene
\begin{equation*}
g (y, \tau ) := \left( \tfrac{\phi^2}{4 \pi} - 4 \phi + 16 \pi \right)\sim 16 \pi \left( \left( \tfrac{y}{\varepsilon ( \tau )} \right)^2 +1 \right)^{-1}
\end{equation*}
cuando $\tau \rightarrow \infty$. Así, obtener la aproximación $g (y, \tau ) \sim 16 \pi^2 \varepsilon ( \tau )^2 \delta (y)$ para $\varepsilon ( \tau ) \ll y \leq 1$, donde $\delta (y)$ es una delta de Dirac centrada en $y = 0$, resulta un paso clave de la demostración. La lógica detrás de esta elección consiste en reemplazar $g (y , \tau )$ en \eqref{ec:Super_ec_psi} por una función de tipo $\gamma \varepsilon ( \tau )^2 \delta (y )$, donde $\gamma$ es elegida de tal forma que $\gamma {\displaystyle \int}_{\mathbb{R}^2} \varepsilon ( \tau )^2 \delta (y) \, dy = {\displaystyle \int}_{\mathbb{R}^2} g (y, \tau ) \, dy$, lo cual implica $\gamma = 16 \pi^2$. 

De aquí en adelante, deberemos asumir $g (y, \tau ) \gg e^{- \tau } \phi (y, \tau )$ para $y \leq 1$ y $\tau \gg 1$, i.e. $\varepsilon ( \tau )^2 \gg e^{- \tau}$ para $\tau \gg 1$. De esta manera, esperamos:
\begin{supo}
La función $\psi (y, \tau )$ se comporta asintóticamente como la solución de
\begin{equation}\label{ec:Super_Sup_1}
\psi_{\tau} = \psi_{yy} + \left( \frac{1}{y} - \frac{y}{2} \right) \psi_y + \psi + \gamma \varepsilon ( \tau )^2 \delta (y) =: A \psi + \gamma \varepsilon ( \tau )^2 \delta (y)
\end{equation}
para $\varepsilon (\tau ) \ll y \leq 1$ y $\tau \gg 1$, donde $\gamma = 16 \pi^2$.
\end{supo}
Para analizar \eqref{ec:Super_Sup_1} añadimos un par de apuntes sobre el espectro de $A$. Tomemos,
\begin{equation*}
L^2_{ \omega , r } ( \mathbb{R}^+ ) = \left\lbrace f \in L^2_{\mathrm{loc}} (\mathbb{R}^+) : \| f \|^2 = \int_0^{\infty} y | f (y) |^2 e^{-y^2 /4} \, dy < \infty \right\rbrace, 
\end{equation*}
que es un espacio de Hilbert con la norma $\| f \|^2 = \left\langle f, f \right\rangle = {\displaystyle \int}_0^{\infty} y | f(y) |^2 e^{-y^2 /4} \, dy < \infty$. 

Para cualquier entero positivo $k$, uno define los espacios de Hilbert $H^k_{\omega, r } ( \mathbb{R}^+ )$ de manera directa. Entonces, resulta que el operador radial $A$ en \eqref{ec:Super_Sup_1} es auto-adjunto en $L^2_{\omega , r} ( \mathbb{R}^+ )$. Además, los autovalores de $A$ son los valores de la sucesión $ \lambda_k = 1 - k$, $k=0, 1, 2, \cdots$, y sus correspondientes autofunciones pueden escribirse de la forma $\varphi_k (y) = c_k L_k \left( \frac{y^2}{4} \right) = \left( \frac{1}{4 \pi k!} \right)^{1/2} L_k \left( \frac{y^2}{4} \right)$ donde $L_k (\xi )$ denota al $k$-ésimo polinomio de Laguerre, con $c_k$ una constante de normalización elegida para tener $\| \varphi_k \| = 1$. Con esto en mente, ahora escribimos la solución $\psi (y, \tau )$ de \eqref{ec:Super_Sup_1} de la forma
\begin{equation}\label{ec:Super_base_autofunc}
\psi (y, \tau ) = a_0 ( \tau ) \varphi_0 (y) + a_1 (\tau ) \varphi_1 (y ) + Q (y, \tau ),
\end{equation}
así, separamos en el análisis la parte generada por las autofunciones de autovalores negativos del resto. La autofunción $\varphi_0$ es la única con autovalor positivo, por ello, para construir nuestro ejemplo evitando problemas impondremos en la Suposición \ref{sup:Super_sup_2} que $a_0 ( \tau )$ tienda a $0$ cuando $ \tau \rightarrow \infty$ y que $a_1 ( \tau ) \varphi_1 ( \tau )$ sea el término principal. Los dos primeros coeficientes de Fourier de \eqref{ec:Super_base_autofunc} satisfacen
\begin{equation}\label{ec:Super_a0_a1}
\left\lbrace 
\begin{array}{ll}
\dot{a_0} & = a_0 + \gamma \varepsilon (\tau )^2 \left\langle  \varphi_0 , \delta (y ) \right\rangle  \\
\dot{a_1} & = \gamma \varepsilon (\tau )^2 \left\langle  \varphi_1, \delta (y) \right\rangle ,
\end{array} \right. 
\end{equation}
mientras que el término $Q (y, \tau)$ es tal que 
\begin{equation}
Q_{\tau} = A Q + \gamma \varepsilon ( \tau )^2 \left( \delta (y) - \sum_{k=0}^1 \varphi_k \left\langle \varphi_k , \delta (y) \right\rangle  \right)
\end{equation}
y $\left\langle  Q, \varphi_k \right\rangle  = 0$ para $k=0, 1$.

De aquí en adelante, utilizaremos $\gamma$ para denotar la constante $16 \pi^2$. También añadimos una segunda suposición.
\begin{supo}\label{sup:Super_sup_2}
El término principal en \eqref{ec:Super_base_autofunc} es $a_1 ( \tau ) \varphi_1 (y )$, así que la evolución en tiempo de $\psi (y, \tau )$ está dirigida por la autofunción correspondiente al autovalor cero. Además, esperamos 
\begin{equation}
\begin{array}{ll}
| \dot{\varepsilon} ( \tau ) | \ll \varepsilon ( \tau ) & \, \, \, \, \mathrm{cuando} \, \, \, \, \tau \rightarrow \infty , \\
Q(y , \tau ) \sim \gamma \varepsilon ( \tau )^2 F(y) & \, \, \, \, \mathrm{cuando} \, \, \, \, \tau \rightarrow \infty ,
\end{array}
\end{equation}
para alguna función $F(y)$. Entonces, resulta que $F(y)$ satisface
\begin{equation}\label{ec:Super_F}
A F + \left( \delta (y) - \sum_{k=0}^1 \varphi_k < \varphi_k , \delta (y) > \right) = 0, \, \, \, \mathrm{con} \, \, \,  \langle F, \varphi_k \rangle = 0  \, \, \, \mathrm{para} \, \, \,  k=0, 1.
\end{equation}
\end{supo}
Integrando por partes \eqref{ec:Super_F} obtenemos,
\begin{equation}
F(y) = - \frac{1}{2 \pi} \log y + B + \mathcal{O} (y^2 | \log y | ) \, \, \, \mathrm{cuando} \, \, \, y \rightarrow 0,
\end{equation}
donde la constante $B$ se determina gracias a las condiciones de ortogonalidad dadas en \eqref{ec:Super_F}. Dado que esperamos $a_k (\tau ) \rightarrow 0$ para $k=0, 1$ cuando $\tau \rightarrow \infty$ en regiones donde $y = \mathcal{O} (1)$, deducimos de \eqref{ec:Super_a0_a1} que
\begin{equation}\label{ec:Super_Expresion_a0_a1}
a_k ( \tau ) \sim - \gamma \varphi_k (0) \int_{\tau }^{\infty} \varepsilon (s)^2 e^{(1-k)( \tau - s )} \, ds, \, \, \, \mathrm{para} \, \, \, k=0,1 \, \, \, \mathrm{cuando} \, \, \, \tau \rightarrow \infty.
\end{equation}
Con todo esto, estamos preparados para unir las expansiones interiores y exteriores de $W( y, \tau )$. Para comenzar, observamos que si $y \sim \varepsilon ( \tau )$ y $ \tau \gg 1$,
\begin{equation}\label{ec:Super_inner_expansion_W}
\begin{array}{ll}
W (y, \tau ) & \sim 8 \pi {\displaystyle \int}_0^y r \left( \frac{r}{\varepsilon (\tau )} \right)^2 \left( \left( \frac{r}{\varepsilon (\tau )} \right)^2 +1 \right)^{-1} \, dr \\
  & = 4 \pi y^2 - 4 \pi \varepsilon^2 \log \left( 1 + \left( \frac{y}{\varepsilon (\tau )} \right)^2 \right) .
\end{array}
\end{equation}
Por otro lado, de \eqref{ec:Super_base_autofunc}-\eqref{ec:Super_Expresion_a0_a1} obtenemos la siguiente expansión exterior para $W$:
\begin{equation}\label{ec:Super_outer_expansion_W}
\begin{array}{ll}
W (y, \tau ) \sim & 4 \pi y^2 - \gamma \sum_{k=0}^1 \varphi_k (0) \varphi_k (y) {\displaystyle \int}_{\tau }^{\infty } \varepsilon (s)^2 e^{(1-k)( \tau -s )} \, ds \\
& + \gamma \varepsilon (\tau )^2 \left( - \frac{1}{2 \pi} \log y + B + \cdots \right),
\end{array}
\end{equation}
válido para $\varepsilon ( \tau ) \ll y < 1$ y $ \tau \gg 1 $. Uniendo \eqref{ec:Super_inner_expansion_W} y \eqref{ec:Super_outer_expansion_W} en las regiones $y = \varepsilon ( \tau )^{1/2}$ se deduce la siguiente ecuación integral para $\varepsilon ( \tau )$:
\begin{equation}\label{ec:Super_glueing}
\begin{array}{ll}
4 \pi \varepsilon ( \tau )^2 \log ( \varepsilon ( \tau )^2 ) = & - \gamma \varphi_0 ( 0)^2 {\displaystyle \int}_{\tau}^{\infty} e^{ \tau -s} \varepsilon (s )^2 \, ds \\
 & - \gamma \varphi_1 (0)^2 {\displaystyle \int}_{ \tau}^{\infty} \varepsilon (s)^2 \, ds + B \gamma \varepsilon ( \tau )^2.
\end{array}
\end{equation}
Un análisis de la ecuación \eqref{ec:Super_glueing} nos da como resultado finalmente a
\begin{equation}\label{ec:Super_estimacion_epsilon}
\varepsilon ( \tau )^2 \sim K e^{- \sqrt{2} \tau^{1/2}} (\tau )^{\frac{1}{2 \sqrt{\tau}} - \frac{1}{2} } \, \, \, \mathrm{cuando} \, \, \, \tau \rightarrow \infty ,
\end{equation}
para alguna $K = K(B) > 0$, y de donde se siguen las estimaciones de los Teoremas \ref{th:Herrero_Velazquez_1} y \ref{th:Herrero_Velazquez_2}. De esta forma, \eqref{ec:Super_estimacion_epsilon} proporciona la expansión para $R(t)$ en \eqref{ec:HVT2_perfil_R}. Las expresiones \eqref{ec:Super_T1_rho}, \eqref{ec:Super_T1_f} y \eqref{ec:HV2_forma_rho} se siguen de las expansiones exteriores e interiores \eqref{ec:Super_outer_expansion_W} y \eqref{ec:Super_inner_expansion_W}.

\chapter{Técnicas de entropía}\label{ch:entropia}

A lo largo de estas notas se ha recurrido a herramientas relacionadas con el transporte óptimo, la distancia de Wasserstein y el estudio de entropías. Las técnicas que hemos desarrollado no son exclusivas para la ecuación de Keller-Segel y pueden ser empleadas en otras similares a \eqref{ec:KS_clasico_una_linea}. Este Capítulo tiene como propósito profundizar más en ellas y nosotros nos centraremos en la ecuación de Fokker-Planck
\begin{equation}\label{ec:FP}
\dfrac{\partial \rho}{\partial t} = \nabla \cdot \left( D ( \nabla \rho + \rho \nabla V) \right), \, \, \, \, \, \, \, \, \, \, \, \, \, \, \, \, \, \, \, \, t \geq 0,\, \, \, \, \, \, \, \, x \in \mathbb{R}^n ,
\end{equation}
donde la difusión viene dada por $D =D(x)$ una matriz simétrica, localmente definida positiva uniformemente, y $V(x)$ un potencial limitante (\textit{confining}). Con más generalidad, analizaremos la ecuación,
\begin{equation}\label{ec:Car_McC_Vil}
\dfrac{\partial \rho}{\partial t} = \Delta \rho + \nabla \cdot \left( \rho  \nabla (  \rho \ast W) \right), \, \, \, \, \, \, \, \, \, \, \, \, \, \, \, \, \, \, \, \, t \geq 0,\, \, \, \, \, \, \, \, x \in \mathbb{R}^n ,
\end{equation}
donde $W$ es un potencial de interacción simétrico $C^2$ que satisface
\begin{equation*}
D^2 W (z) \geq K |z|^{\gamma} \, \, \, \, \, \, \, \, \, \, \, \mathrm{y} \, \, \, \, \, \, \, \, \, \, \, | \nabla W(z) | \leq C (1+|z|)^{\beta}
\end{equation*}
para algún $\gamma >0$ y algún $\beta \geq 0$ (por ejemplo $W(z) = |z|^3/3$).

Asumimos que existe un $\rho_{\infty}$ que es minimizante de alguna energía asociada a la ecuación que estemos estudiando. Nosotros queremos desarrollar un método que demuestre convergencia exponencial a este minimizante. Uno que se ha mostrado especialmente eficaz es el estudio de la comparación de la \textit{disipación de la entropía} con la \textit{disipación de la disipación de la entropía}, que ha servido para probar resultados sobre la velocidad de convergencia en distintas ecuaciones, la ecuación de Boltzmann \cite{CAR_CAR, CAR_CAR_2}, Fokker-Planck \cite{ARN_MAR_TOS_UNT, CAR_TOS_2}, la ecuación de medios porosos \cite{CAR_TOS} o en varios sistemas de ecuaciones parabólicas degeneradas \cite{CAR_JUN_MAR_TOS_UNT}, entre otros ejemplos. Por la necesidad de utilizar este minimizador único no desarrollamos estas ideas durante el Capítulo \ref{ch:Subcritica}, ya que para la energía $E_{KS}$ que estábamos estudiando solo sabíamos que esta tenía algún mínimo local.

Comenzamos estudiando la ecuación de Fokker-Planck \eqref{ec:FP} porque en ella es más sencillo comprender la esencia del método. En las Secciones \ref{sec:Disipacion_entropia} y \ref{sec:Des_log_Sob} se discuten los resultados obtenidos por Arnold, Markowich, Toscani y Unterreiter en los papers clásicos \cite{ARN_MAR_TOS_UNT_2, ARN_MAR_TOS_UNT} para las ecuaciones de Fokker-Planck \eqref{ec:FP}\footnote{en \cite{RIS} se puede aprender sobre esta ecuación}, donde se demuestra velocidad de convergencia exponencial al estado estacionario, usando para ello las desigualdades logarítmicas de Sobolev y el Teorema de Bakry-Emery \cite{BAK_EME}. En la Sección \ref{sec:Displacement_convexity} abordamos el mismo problema desde otro punto de vista para conseguir una formulación equivalente. Para ello definimos el concepto de \textit{displacement convexity}, introducido por McCann en \cite{MCC} y desarrollado por Otto y Villani en \cite{OTT, OTT_VIL}. En estas tres secciones seguimos las notas \cite{MAR_VIL}. 

Tras haber estudiado estas técnicas para la ecuación de Fokker-Planck las aplicaremos a la ecuación \eqref{ec:Car_McC_Vil} en la Sección \ref{sec:Car_McC_Vil_Convergencia_exponencial}, siguiendo el trabajo de Carrillo, McCann y Villani \cite{CAR_MCC_VIL}.

Tomemos la versión simplificada de la ecuación de Fokker-Planck,
\begin{equation}\label{ec:FP_simplificado}
\dfrac{\partial \rho}{\partial t} = \nabla \cdot \left(  \nabla \rho + \rho \nabla V \right), \, \, \, \, \, \, \, \, \, \, \, \, \, \, \, \, \, \, \, \, t \geq 0,\, \, \, \, \, \, \, \, x \in \mathbb{R}^n \, \, \, \mathrm{o} \, \, \, M,
\end{equation}
con $M$ una variedad Riemanniana y $V \in C^2$ (elegido de esta forma para evitar problemas sobre regularidad). Estudiamos el problema de Cauchy con dato inicial
\begin{equation*}
\rho (t=0, \cdot ) = \rho_0; \, \, \, \, \, \, \, \, \, \, \, \, \, \, \, \, \rho_0 \geq 0, \, \, \, \, \, \, \, \, \int \rho_0=1.
\end{equation*}
Observamos que $\rho_{\infty} = e^{-V}$ es un estado estacionario de \eqref{ec:FP_simplificado} (añadiendo una constante si fuera necesario) y puede asumirse que $e^{-V}$ es una distribución de probabilidad. Con un cambio de variables $\rho = h e^{-V}$ se obtiene una formulación equivalente de \eqref{ec:FP_simplificado}
\begin{equation}\label{ec:FP_simplificado_h}
\dfrac{\partial h}{\partial t} = \Delta h - \nabla V \cdot \nabla h, \, \, \, \, \, \, \, \, \, \, \, \, \, \, \, \, \, \, \, \, t \geq 0,\, \, \, \, \, \, \, \, x \in \mathbb{R}^n \, \, \, \mathrm{o} \, \, \, M.
\end{equation}
Se espera que la solución del problema de Cauchy converja al estado de equilibrio $e^{-V}$ y nos gustaría estimar la velocidad de convergencia en términos del dato inicial. Con este fin enunciamos el Teorema \ref{th:Velocidad_convergencia_FP}
\begin{teor}\label{th:Velocidad_convergencia_FP}
Sea $\rho$ una densidad que resuelve \eqref{ec:FP_simplificado} para dato inicial $\rho_0$, se puede alcanzar la estimación
\begin{equation}
\rho_0 \in L^2 \left( e^V \right) \, \, \Longrightarrow \, \, \| \rho (t, \cdot ) - e^{-V} \|_{L^2 \left(e^V \right)} \leq e^{- \lambda t} \| \rho_0 - e^{-V} \|_{L^2\left( e^V \right)} .
\end{equation}
\end{teor}

En la Sección \ref{sec:Disipacion_entropia} introducimos más a fondo las técnicas y las entropías que usaremos.

\section{Disipación de la entropía}\label{sec:Disipacion_entropia}

En lugar de investigar el decaimiento de $h$ en norma $L^2 \left( e^{-V} \right)$, podríamos considerar una variedad de funcionales controlando la distancia entre $h$ y $1$. Para cualquier función convexa $\phi$ en $\mathbb{R}$, uno puede comprobar que 
\begin{equation}
\int \phi (h) e^{-V} \, dx = \int \phi \left( \dfrac{\rho}{e^{-V}} \right) e^{-V} \, dx
\end{equation}
define un funcional de Lyapunov para \eqref{ec:FP_simplificado_h}, o equivalentemente \eqref{ec:FP_simplificado}, ya que,
\begin{equation}\label{ec:FP_Lyapunov}
\dfrac{d}{dt} \int \! \! \! \int \phi (h) e^{-V} \, dx = - \int \phi'' (h) | \nabla h |^2 e^{-V} \, dx \leq 0.
\end{equation}
Es interesante elegir $\phi (h) = h \log h - h +1$. Hemos añadido el término $-h+1$, que satisface $\int (h-1) e^{-V} = 0$, para conseguir $\phi (h) \geq 0$. Encontramos 
\begin{equation}\label{ec:FP_entropia_relativa}
\int \phi (h) e^{-V} = \int \rho \log \dfrac{\rho}{e^{-V}} = \int \rho ( \log \rho + V ) =: H \left( \rho \, | \, e^{-V} \right),
\end{equation} 
un funcional que es bien conocido. En teoría cinética \cite{KAC} se llama la \textit{energía libre}, mientras que en teoría de la información \cite{CSI_2, KUL, PER_2} se conoce como la \textit{entropía relativa} (Kullback) de $\rho$ con respecto a $e^{-V}$ (i.e., la medida $\rho \, dx$ con respecto a la medida $e^{-V} \, dx$). \eqref{ec:FP_entropia_relativa} es una \textit{entropía relativa} que denotamos por $H \left( \rho \, | \, e^{-V} \right)$ (utilizamos la notación estándar de la entropía de Boltzmann).

La entropía relativa es un buen candidato para controlar la distancia entre dos distribuciones de probabilidad gracias a la desigualdad\footnote{\eqref{ec:Csiszar_Kullback_Pinsker} se conoce en análisis como la desigualdad Csiszár-Kullback \cite{CSI, KUL}, y como la desigualdad de Pinsker \cite{PIN_2} en probabilidad (cf. Apéndice A de \cite{REI_WIL} para ver la historia de esta desigualdad).} elemental
\begin{equation}\label{ec:Csiszar_Kullback_Pinsker}
H \left( \rho \, | \, \tilde{\rho} \right) \geq \dfrac{1}{2} \| \rho - \tilde{\rho} \|^2_{L^1} .
\end{equation}
Por \eqref{ec:FP_Lyapunov}, si $\rho$ es solución de \eqref{ec:FP_simplificado}, entonces
\begin{equation}\label{ec:FP_funcional_Fisher}
\dfrac{d}{dt} H \left( \rho \, | \, e^{-V} \right) = - \int \rho \left| \nabla \left( \log \dfrac{\rho}{e^{-V}} \right) \right|^2 \, dx =: - I \left( \rho \, | \, e^{-V} \right) ,
\end{equation}
con $I$ denominada la \textit{disipación de la entropía} o información de Fisher. Los cálculos \eqref{ec:FP_funcional_Fisher} nos ayudan a investigar e identificar la tendencia al equilibrio para \eqref{ec:FP_simplificado}.

\section{Desigualdades logarítmicas de Sobolev}\label{sec:Des_log_Sob}

El último paso que necesitamos para probar convergencia exponencial puede darse gracias a las desigualdades logarítmicas de Sobolev, y es por eso qué necesitamos entender cómo se aprovechan para nuestro problema y bajo que hipótesis podemos usarlas.

Sea $\gamma_{\sigma} (x) = (2 \pi \sigma)^{-n/2} e^{- |x|^2 / (2 \sigma)}$, la Gaussiana centrada con varianza $\sigma$. La desigualdad logarítmica de Stam-Gross de Sobolev \cite{GRO, STA} asegura que para cualquier distribución de probabilidad $\rho$ (absolutamente continua con respecto a $\gamma_{\sigma}$),
\begin{equation}\label{ec:FP_des_log_Sob}
H \left( \rho \, | \, \gamma_{\sigma} \right) \leq \dfrac{\sigma}{2} I \left( \rho \, | \, \gamma_{\sigma} \right) .
\end{equation}
Esta se llama desigualdad logarítmica de Sobolev porque puede reescribirse como
\begin{equation*}
\int u^2 \log u^2 \, d \gamma - \left( \int u^2 \, d \gamma \right) \log \left( \int u^2 \, d \gamma \right) \leq 2 \int | \nabla u |^2 \, d \gamma .
\end{equation*}
Con \eqref{ec:FP_des_log_Sob} continuamos con el estudio de la tendencia al equilibrio para \eqref{ec:FP_simplificado}. 

Comenzamos con un caso sencillo, $V = \frac{|x|^2}{2 \sigma}$ un potencial cuadrático y más adelante intentaremos generalizar el argumento. Por \eqref{ec:FP_funcional_Fisher}, si $\rho$ es solución de 
\begin{equation}\label{ec:FP_Caso_particular}
\dfrac{\partial \rho}{\partial t} = \nabla \cdot \left( \nabla \rho + \rho \dfrac{x}{\sigma} \right),
\end{equation}
entonces $\rho$ satisface una estimación de \textit{decaimiento exponencial de la entropía relativa},
\begin{equation}\label{ec:FP_entrop_relativa_decaimiento_exponencial}
H \left( \rho (t, \cdot ) \, | \, \gamma_{\sigma} \right) \leq H \left( \rho_0 \, | \, \gamma_{\sigma} \right) e^{- 2t / \sigma} .
\end{equation}
La solución de \eqref{ec:FP_entropia_relativa} se puede calcular, lo que motiva intentar entender cómo se pueden extender estos resultados a condiciones más generales.

Por definición decimos  que una medida de probabilidad $e^{-V}$ satisface la desigualdad logarítmica de Sobolev con constante $\lambda > 0$ si para toda medida de probabilidad $\rho$,
\begin{equation}
H \left( \rho \, | \, e^{-V} \right) \leq \dfrac{1}{2 \lambda} I \left( \rho \, | \, e^{-V} \right).
\end{equation}
Con cálculos análogos a los que acabamos de desarrollar nos damos cuenta de que si $e^{-V}$ satisface la desigualdad logarítmica de Sobolev con constante $\lambda$, entonces la solución de Fokker-Planck (con $V$ el potencial de confinamiento) se va al equilibrio en entropía relativa, con una velocidad de al menos $e^{-2 \lambda t}$. 

Tras esto surge una pregunta natural, \textit{¿cuáles son las medidas de probabilidad que satisfacen las desigualdades logarítmicas de Sobolev?} En 1985, Bakry y Emery \cite{BAK_EME} le dieron una respuesta enunciando el Teorema.
\begin{teor}[Bakry-Emery \cite{BAK_EME}]\label{th:Bakry_Emery}
Sea $e^{-V}$ una medida de probabilidad en $\mathbb{R}^n$ (resp. una variedad Riemanniana $M$), tal que $D^2 V \geq \lambda I_n$ (resp. $D^2 V + \mathrm{Ric} \geq \lambda I_n$), donde $I_n$ es la matriz identidad de dimensión $n$ (y $\mathrm{Ric}$ el tensor de curvatura de Ricci para la variedad $M$). Entonces, $e^{-V}$ satisface una desigualdad logarítmica de Sobolev con constante $\lambda$.
\end{teor}
(En el caso Riemanniano, $D^2 V$ hace referencia al Hessiano de $V$). Además, todavía hay espacio para la perturbación en este Teorema, como puede verse del lema de perturbación estándar de Holley-Stroock \cite{HOL_STR}: si $V$ es de la forma $V_0 + v$, donde $v \in L^{\infty}$ y $e^{-V_0}$ satisface una desigualdad logarítmica de Sobolev con constante $\lambda$, entonces, $e^{-V}$ también satisface una desigualdad logarítmica de Sobolev, con constante $\lambda e^{- \mathrm{osc} (v)}$, con $\mathrm{osc} (v) = \sup v - \inf v$. Al combinar el Teorema de Bakry-Emery con el lema de Holley-Strock, se pueden generar numerosas medidas de probabilidad que satisfacen la desigualdad logarítmica de Sobolev (para enunciados más generales con una matriz de difusión $D$ no-constante nos referimos a \cite{ARN_MAR_TOS_UNT_2}). En estas notas no daremos una demostración del Teorema de Bakry-Emery, pero puede seguirse el trabajo de Arnold, Markowich, Toscani y Unterreiter \cite{ARN_MAR_TOS_UNT_2} para una demostración desde un punto de vista físico.

\section{Displacement convexity}\label{sec:Displacement_convexity}

El Teorema \ref{th:Bakry_Emery} nos ha permitido dar una condición suficiente sobre decaimiento exponencial de la entropía relativa, si $D^2 V \geq \lambda I_n$ entonces hay decaimiento exponencial. En este capítulo, queremos encontrar otra forma equivalente de plantear esta condición que sea más sencilla de comprobar cuando quieres aplicar el resultado a casos concretos. Nuestro objetivo es expresar esta condición suficiente solo en términos de de la entropía relativa $H$. Esto se puede hacer gracias a la \textit{displacement convexity} (o \textit{convexidad de Wasserstein}), herramienta introducida por McCann en \cite{MCC} y desarrollada por Otto y Villani en \cite{OTT, OTT_VIL}. 

La displacement convexity es un concepto muy interesante en si mismo. El análisis de funcionales de energía tiene un papel importante en varios campos de las matemáticas (en este trabajo, por ejemplo, aparecen varias veces). Normalmente, intentamos determinar la existencia de estados estacionarios, particularmente optimizadores y sus propiedades: unicidad, estabilidad, simetría... Disponer de alguna noción de convexidad ayuda a resolver estos problemas. La displacement convexity es una nueva noción de convexidad que surge de manera natural de la interpolación de medidas de probabilidad en $\mathbb{R}^d$.

Presentadas las ventajas de la displacement convexity y sus razones para estudiarlas, definimos lo que es. Para ello, recordamos\footnote{En la Sección \ref{sec:Sub_JKO} hemos normalizado el paso en tiempo para que fuese $\tau$ y no $1$} como interpolamos dos medidas de probabilidad (hecho en la Sección \ref{sec:Sub_JKO}). Sean $\rho_0$ y $\rho_1$ dos medidas de probabilidad, existe $\varphi$ convexo en $\mathbb{R}^d$ tal que $\nabla \varphi \, \sharp \, \rho_0 = \rho_1$. Sea $id$ la aplicación identidad en $\mathbb{R}^d$. Para tiempo $ t \in [0, 1]$, el interpolante $\rho_t \in \mathcal{P} \left( \mathbb{R}^d \right)$ entre $\rho_0$ y $\rho_1$ queda definido por
\begin{equation*}\label{eq:Interpolacion_Wasserstein}
\rho_t := [ (1-t) id + t \nabla \varphi ] \, \sharp \, \rho_0.
\end{equation*}
Esta interpolación recibe el nombre de \textit{displacement interpolation}.
 
Entonces, decimos que un funcional cualquiera $J$ es displacement convex si $t \mapsto J( \rho_t)$ es convexo en $[0,1]$ para todas las parejas de distribuciones de probabilidad $\rho_0$, $\rho_1$. También decimos que $J$ es uniformemente displacement convex con constante $\lambda$ si
\begin{equation*}
\dfrac{d^2}{dt^2} J (\rho_t) \geq \lambda \mathcal{W}_2 (\rho_0, \rho_1 )^2, \, \, \, \, \, \, \, \, \, \, \, \, \, \, \, 0 < t <1.
\end{equation*}
Aplicada a nuestro funcional, si $V$ es uniformemente convexo con constante $\lambda$ tenemos que la entropía relativa $H \left( \rho \, | \, e^{-V} \right) = \int \rho \log \rho + \int \rho V$ define un funcional (uniformemente) displacement convex.

Tal y como ocurre con la convexidad estándar, la noción de displacement convexity puede formularse con argumentos de diferenciabilidad. Siguiendo esta línea de razonamiento, Otto en \cite{OTT} introduce una estructura en el conjunto de medidas de probabilidad, de tal forma que la displacement convexity (resp. uniformemente displacement convexity) del funcional $J$ es equivalente a la no-negatividad (resp. positividad uniforme) del Hessiano de $J$.

Según hemos visto en el apéndice \ref{ap:Transporte_optimo}, dados un funcional $E$ en el conjunto de medidas de probabilidad, los cálculos formales de Otto permiten definir de forma natural un funcional $| \mathrm{grad} \, E |^2$ por
\begin{equation*}
| \mathrm{grad} \, E |^2 (\rho) = \int \rho | \nabla \psi |^2,
\end{equation*}
donde $\psi [= \delta E / \delta \rho]$ es la solución de [$'$ derivada funcional]
\begin{equation*}
\int \psi g = E' (\rho) \cdot g.
\end{equation*}
Para el caso particular que nos atañe, $E (\rho) = \int \rho \log \rho + \int \rho V$, i.e. $E (\rho) = H \left( \rho \, | \, e^{-V} \right)$. Uno comprueba fácilmente que $| \mathrm{grad} E |^2 ( \rho) = I \left( \rho \, | \, e^{-V} \right)$. Suponiendo ahora que $D^2 V \geq \lambda I_n$ el resultado de Bakry-Emery puede reescribirse como
\begin{equation*}
| \mathrm{grad} \, H \left( \cdot \, | \, e^{-V} \right) |^2 \geq 2 \lambda H \left( \cdot \, | \, e^{-V} \right).
\end{equation*}  
Esta desigualdad, formalmente, es una consecuencia muy sencilla del hecho de que la entropía relativa es uniformemente displacement convex con constante $\lambda$. Pese a su sencillez, nos permite reescribir el Teorema de Bakry-Emery \ref{th:Bakry_Emery} con una condición que solo depende de la entropía relativa $H$, tal y como queríamos.

\section{Convergencia exponencial en $L^1$ para la ecuación \eqref{ec:Car_McC_Vil}}\label{sec:Car_McC_Vil_Convergencia_exponencial}

En esta última Sección aprovecharemos toda la teoría sobre la comparación de la \textit{disipación de la entropía} con la \textit{disipación de la disipación de la entropía} desarrollada en las Secciones anteriores, para probar el siguiente resultado.

\begin{teor}[Carrillo-McCann-Villani \cite{CAR_MCC_VIL}]
Sea $\rho = (\rho_t )_{t \geq 0}$ solución de 
\begin{equation*}
\dfrac{\partial \rho}{\partial t} = \Delta \rho + \nabla \cdot \left( \rho  \nabla (  \rho \ast W) \right), \, \, \, \, \, \, \, \, \, \, \, \, \, \, \, \, \, \, \, \, t \geq 0,\, \, \, \, \, \, \, \, x \in \mathbb{R}^n ,
\end{equation*}
donde $W$ es un potencial de interacción simétrico $C^2$ satisfaciendo
\begin{equation*}
D^2 W (z) \geq K |z|^{\gamma} \, \, \, \, \, \, \, \, \, \, \, \mathrm{y} \, \, \, \, \, \, \, \, \, \, \, | \nabla W(z) | \leq C (1+|z|)^{\beta}
\end{equation*}
para algún $\gamma >0$ y algún $\beta \geq 0$. Sea 
\begin{equation*}
F [\rho] = \int \rho \log \rho \, dx + \dfrac{1}{2} \int W (x-y) \, d \rho (x) d \rho (y),
\end{equation*}
y sea $\rho_{\infty}$ el minimizador único de $F$ con el mismo centro de masas que $\rho$. Entonces, para cada $t_0 >0$ existen constantes $C_0, \lambda_0 > 0$ tales que
\begin{equation*}
t \geq t_0 \, \Rightarrow \, \| \rho_t - \rho_{\infty} \|_{L^1} \leq C_0 e^{-\lambda_0 t }.
\end{equation*}
\end{teor}
Antes de comenzar con la prueba, queremos señalar que Carrillo, McCann y Villani prueban la unicidad del minimizador de $F$ en el trabajo \cite{CAR_MCC_VIL}, utilizando los resultados de McCann \cite{MCC}. Nosotros no discutiremos sobre esto aquí. Por lo demás, seguiremos la misma prueba del Teorema dada por Carrillo, McCann y Villani en \cite{CAR_MCC_VIL}.

Por la displacement convexity de $F$ se sigue como en \cite{OTT_VIL_2},
\begin{equation*}
F [\rho_t] \leq \frac{\mathcal{W}_2 (\rho_0, \rho_{\infty})^2}{4 t} \, \, \, \, \, \, \, \, \forall t >0.
\end{equation*}
Por definición, $\mathcal{W}_2 (\rho_0, \rho_{\infty})^2 \leq 2 \left( \int |x|^2 \, d\rho_0 (x) + \int |x|^2 \, d\rho_{\infty} (x) \right)$, nos da una cota \textit{a priori} en $F [\rho_t ]$ para $t > t_0 \geq 0$ que solo depende de $\int |x|^2 \, d \rho_0 (x)$. Combinando este resultado con el Teorema \ref{th:Carrillo_McCann_Villani_2.4} podemos recuperar convergencia exponencial de $F [ \rho_t ]$ a $F [\rho_{\infty}]$, y también 
\begin{equation}\label{ec:Subcritica_exp_Primera_desigualdad}
\mathcal{W}_2 (\rho_t , \rho_{\infty} ) \leq C_1 e^{-\lambda_1 t}, \, \, \, \, \, \, \, \, \, \, \, \, \, \, \, \, \, \, t \geq t_0 .
\end{equation}

Adaptado al caso particular que estamos estudiando, el Teorema 2.4 de \cite{CAR_MCC_VIL} enuncia:
\begin{teor}\label{th:Carrillo_McCann_Villani_2.4}
La ecuación \eqref{ec:Car_McC_Vil} y $F$ satisfacen las propiedades,
\begin{enumerate}[label=\Roman*)]
\item Existe un único minimizador $\rho_{\infty}$ para $E_{KS}$ en la clase de medidas de probabilidad $\rho$ con centro de masas $\theta$. Este minimizador resulta ser el único estado estacionario para \eqref{ec:KS_clasico_una_linea} en la misma clase de medidas de probabilidad.
\item Para todas las densidades de probabilidad $\rho$ tales que $E_{KS} < + \infty$, existe una constante $\lambda = \lambda (\rho )$, que solo depende de $E_{KS}$ y tal que las desigualdades HWI se satisfacen.\label{th:Carrillo_McCann_Villani_2.4_ii}
\item Si $\lambda_0$ es la constante asociada a $\rho_0$ por \ref{th:Carrillo_McCann_Villani_2.4_ii}, existen soluciones $\left( \rho_t \right)_{ t \geq 0}$ de \eqref{ec:KS_clasico_una_linea} tales que $t \mapsto \rho_t$ es continuo (en el sentido de distribución) y que satisfacen una estimación de decaimiento
\begin{equation*}
F [ \rho_t ] - F [\rho_{\infty} ] \leq e^{-2 \lambda_0 t} \left( F [\rho_0 ] - F [\rho_{\infty} ] \right) .
\end{equation*}
\end{enumerate}
\end{teor}

El resto del argumento consiste en transformar esta información de convergencia débil en convergencia fuerte por medio de una interpolación correcta. Para ello introducimos el \textit{funcional de entropía de disipación}
\begin{equation*} 
D (\rho) = \int \nabla | \log \rho - \rho \ast W|^2 \, d \rho (x) . 
\end{equation*}
Tomando otra vez la prueba en \cite{OTT_VIL_2} uno también tiene la estimación
\begin{equation}\label{ec:Subcritica_exp_dis_Was}
D (\rho_t) \leq \frac{\mathcal{W}_2 ( \rho_0, \rho_{\infty})^2}{t^2}
\end{equation}
que nos da una cota uniforme en $D (\rho_t)$ para $t \geq t_0$.

Por otro lado, un análisis de la ecuación de Euler-Lagrange asociado a la minimización de $F$ demuestra que $\rho_{\infty}$ tiene una densidad estrictamente positiva, satisfaciendo 
\begin{equation}\label{ec:Subcritica_exp_forma_compacta}
\log \rho_{\infty} - \rho_{\infty} \ast W = \mu \in \mathbb{R}.
\end{equation}
En particular,
\begin{equation*}
\nabla \log \rho_{\infty} = \nabla ( \rho_{\infty} \ast W ),
\end{equation*}
y
\begin{equation}\label{ec:Subcritica_exp_des_disip}
\begin{array}{lll}
{\displaystyle \int} \left| \nabla \log \left( \frac{\rho_t}{\rho_{\infty}} \right) \right|^2 \, d \rho_t (x) & = {\displaystyle \int} | \nabla \log \rho_t - \nabla ( \rho_{\infty} \ast W )|^2 \, d \rho_t (x) \\
 & \leq 2 {\displaystyle \int} \left[ | \nabla ( \log \rho_t - \rho_t \ast W)|^2 + |(\rho_t + \rho_{\infty} ) \ast \nabla W |^2 \right] \, d \rho_t (x) \\
 & \leq 2 D( \rho_t ) + C \left( {\displaystyle \int} (1 + |x|^{2 \beta} ) \, d \rho_t (x) \right)^2 .
\end{array}
\end{equation}
Una estimación del momento (siguiendo \cite{DES_VIL}) nos lleva a la desigualdad diferencial
\begin{equation*}
\begin{array}{ll}
\dfrac{d}{dt} {\displaystyle \int} (1 + |x|^2)^{\beta} \, d \rho_t (x) & \leq C - K {\displaystyle \int} |x|^{2 \beta + \gamma} \, d \rho_t (x) \\
 & \leq C - K \left[ {\displaystyle \int} (1 + |x|^2)^{\beta} \, d \rho_t (x) \right]^{1 + \frac{\gamma}{2 \beta} } ,
\end{array}
\end{equation*}
con $C$ y $K$ denotando varias constantes que dependen solo de $\rho$ a través de un cota superior en $\int |x|^2 \, d \rho_0 (x)$. Esto implica una cota uniforme en $\int |x|^{2 \beta} \, d \rho_t (x)$ cuando $t \rightarrow \infty$. Combinando esto con \eqref{ec:Subcritica_exp_des_disip} y \eqref{ec:Subcritica_exp_dis_Was} obtenemos la estimación
\begin{equation}\label{ec:Subcritica_exp_penult_desig}
\int \left| \nabla \log \dfrac{\rho_t}{\rho_{\infty}} \right|^2 \, d \rho_t (x) \leq C, \, \, \, \, \, \, \, \, \, \, \, \, \, \, \, \, \, \, t \geq t_0 > 0.
\end{equation}
Ahora, por \eqref{ec:Subcritica_exp_forma_compacta} también deducimos que $\log (d \rho_{\infty} / dx )$ es cóncava (semicóncava sería suficiente para el argumento). Tomamos la desigualdad HWI de \cite{OTT_VIL}, 
\begin{equation*}
\int \log \dfrac{\rho_t}{\rho_{\infty}} \, d \rho_t \leq \mathcal{W}_2 ( \rho_t, \rho_{\infty} ) \sqrt{\int \left| \nabla \log \dfrac{\rho_t}{\rho_{\infty}} \right|^2 \, d \rho_t } ,
\end{equation*}
y recordamos \eqref{ec:Subcritica_exp_Primera_desigualdad} y \eqref{ec:Subcritica_exp_penult_desig}, que implica que el lado izquierdo converge a $0$ de forma exponencial cuando $t \rightarrow \infty$. Para concluir el argumento basta con utilizar la desigualdad
\begin{equation*}
\| \rho_t - \rho_{\infty} \|_{L^1}^2 \leq 2 \int \log \dfrac{\rho_t}{\rho_{\infty}} \, d \rho_t .
\end{equation*}

\appendix
\chapter{Transporte óptimo}\label{ap:Transporte_optimo}

Este Apéndice se escribe con el objetivo de motivar la distancia 2-Wasserstein, que aparece varias veces a lo largo de estas notas. Para ello, relacionamos esta métrica con el problema de transporte óptimo. 

Dividimos el Apéndice en dos Secciones. En la primera, Sección \ref{ap_sec:Monge_Kantorovich}, presentamos el problema del transporte óptimo, con sus dos formulaciones, la de Monge y la de Kantorovich, discutiendo la equivalencia de ambos problemas bajo ciertas circunstancias. En esta Sección seguiremos las notas \cite{AMB_GIG}. 

En la Sección \ref{ap_sec:Transporte_Optimo_Otto}, discutimos sobre las aplicaciones del problema de transporte óptimo dadas por Otto \cite{OTT}, quien dio la intuición para poder tratar este problema con argumentos propios de la geometría. Para esta sección seguimos las notas \cite{VIL_2}.

\section{Formulaciones de Monge y Kantorovich para el problema de transporte óptimo} \label{ap_sec:Monge_Kantorovich}

Sea $\left( X, d \right)$ un espacio métrico completo y separable, denotamos por $\mathcal{P} (X)$ el conjunto de medidas de probabilidad de Borel en $X$.

Si $X$, $Y$ son dos espacios métricos completos y separables, $T: X \rightarrow Y$ una aplicación de Borel, y $\mu \in \mathcal{P} (X)$ una medida, la medida $T \sharp \mu \in \mathcal{P} (Y)$, llamada \textit{push forward de $\mu$ a través de $T$}, se define
\begin{equation*}
T \sharp \mu (E) = \mu (T^{-1} (E)), \, \, \, \, \, \, \, \, \, \, \, \, \, \, \, \, \, \, \, \, \, \, \, \, \forall E \subset Y, \, \, \, \, \mathrm{Borel}.
\end{equation*}
El push forward se caracteriza por el hecho de que 
\begin{equation*}
\int_Y f d T \sharp \mu = \int_Y f \circ T d \mu .
\end{equation*}
para toda función de Borel $f : Y \rightarrow \mathbb{R} \cup \left\lbrace + \infty \right\rbrace$\footnote{Esta identidad debe entenderse con el siguiente sentido: una de las integrales existe (es posible que asignándole el valor $\pm \infty$) si y solo si la otra existe, y en ese caso, ambos valores son iguales.}.

Fijemos una \textit{función de coste} Borel $c: X \times Y \rightarrow \mathbb{R} \cup \left\lbrace + \infty \right\rbrace$. La versión de Monge para el problema de transporte es la siguiente\footnote{La formulación original de G. Monge del problema de transporte \cite{MON} se refería al caso $X=Y=\mathbb{R}^d$ y $c (x, y) = |x-y|$.}: 

\begin{prob}
[Problema de transporte óptimo de Monge, \cite{MON}] \label{prob:Monge} Sea $\mu \in \mathcal{P} (X)$, $\nu \in \mathcal{P} (Y)$. Minimiza
\begin{equation*}
T \mapsto \int_X c \left( x, T(x) \right) d \mu (x)
\end{equation*}
de entre todas los transports maps desde $\mu$ hasta $\nu$, i.e. todos los maps $T$ tales que $T \sharp \mu = \nu$.
\end{prob}

Este Problema se puede resolver, y, para el caso de coste $c$ cuadrático, la respuesta viene dada por el Teorema \ref{th:Benamou_Brenier} de Benamou-Brenier, que veremos con mayor detalle en la Sección \ref{ap_sec:Transporte_Optimo_Otto}.  

No obstante, la formulación del problema de transporte óptimo dada por Monge lleva asociada algunos inconvenientes. Independientemente de la función de coste $c$ elegida, el problema de Monge puede no estar bien definido porque:
\begin{itemize}
\item No exista $T$ admisible (por ejemplo si $\mu$ es una delta de Dirac pero $\nu$ no lo es).
\item La restricción $T \sharp \mu = \nu$ no es cerrada débilmente por sucesiones con respecto a cualquier topología débil razonable.
\end{itemize}
Un ejemplo para el segundo fenómeno es la sucesión $f_n (x) := f(nx)$, donde $f: \mathbb{R} \rightarrow \mathbb{R}$ es $1$-periódica e igual a $1$ en $[0, 1/2)$ y a $-1$ en $[1/2, 1)$, y las medidas $\mu := \mathcal{L}\vert_{[0,1]}$ y $\nu := (\delta_{-1} + \delta_1)/2$. Es inmediato comprobar que $(f_n) \sharp \mu = \nu$ para todo $n \in \mathbb{N}$, $(f_n)$ converge débilmente a la función nula $f \equiv 0$ que satisface $f \sharp \mu = \delta_0 \neq \nu$.

Una forma de evitar esta dificultad es gracias a Kantorovich, que propuso la siguiente forma de relajar el problema:

\begin{prob}
[Problema de transporte óptimo de Kantorovich, \cite{KAN}]\label{prob:Kantorovich} Minimizamos
\begin{equation*}
\gamma \mapsto \int_{X \times Y} c (x, y) d \gamma (x,y)
\end{equation*}
en el conjunto de medidas admisibles $\mathrm{Adm} (\mu, \nu)$ para todos los transport plans $\gamma \in \mathcal{P} (X \times Y)$ desde $\mu$ hasta $\nu$. Aquí el conjunto $\mathrm{Adm} (\mu , \nu)$ se define como el conjunto de medidas de probabilidad de Borel en $X \times Y$ tal que
\begin{equation*}
\gamma (A \times Y ) = \mu (A) \, \, \, \, \, \, \, \forall A \in \mathcal{B} (X), \, \, \, \, \, \, \, \, \, \, \, \, \, \,  \gamma ( X \times B ) = \nu (B) \, \, \, \, \, \, \, \forall B \in \mathcal{B} (Y).
\end{equation*}
Equivalentemente: $\pi^X \sharp \gamma = \mu$, $\pi^Y \sharp \gamma = \nu$, donde $\pi^X$, $\pi^Y$ son las proyecciones naturales de $X \times Y$ en $X$ e $Y$ respectivamente.
\end{prob}

Para entender los transport plans, podemos pensar en ellos como transport maps ``multievaluados'': $\gamma = \int \gamma_x d \mu (x)$, con $\gamma_x \in \mathcal{P} \left( \left\lbrace x \right\rbrace \times Y \right)$. Otra forma sería observar que para $\gamma \in \mathrm{Adm} ( \mu, \nu)$, el valor de $\gamma ( A \times B)$ es la cantidad de masa que se encontraba inicialmente en $A$ y se ha desplazado a $B$.

Existen varias ventajas en la formulación del problema de transporte dada por Kantorovich:
\begin{itemize}
\item $\mathrm{Adm} (\mu, \nu)$ es siempre no vacío (contiene a $\mu \times \nu$).
\item El mínimo siempre existe bajo hipótesis de suavidad para $c$, es suficiente con que $c$ sea lower semicontinuous y acotado por abajo.
\item Transport plans ``incluyen'' transport maps, ya que $T \sharp \mu = \nu$ implica que $\gamma := (I\!d \times T) \sharp \mu$ pertenece a $\mathrm{Adm} (\mu, \nu)$.
\end{itemize}

Por último, estudiamos la relación entre el Problema \ref{prob:Kantorovich} y \ref{prob:Monge}. Es posible demostrar que si $c$ es continuo y $\mu$ no es atómico, entonces
\begin{equation*}
\inf (\mathrm{Monge}) = \min (\mathrm{Kantorovich}),
\end{equation*}
lo que significa que transportar con plans no puede ser estrictamente más barato que transportar con maps. 

Esta igualdad fue demostrada en casos particulares por Gangbo (Apéndice A de \cite{GAN}) y por Ambrosio  (Teorema 2.1 de \cite{AMB}). La demostración del caso generalizado fue obtenida más tarde por Pratelli en \cite{PRA}. Nosotros no daremos la demostración de este resultado en estas notas.

\section{Punto de vista de Otto} \label{ap_sec:Transporte_Optimo_Otto}

En esta parte del Apéndice vamos a presentar la construcción de Otto \cite{OTT}. Esta ofrece una re-interpretación formal del problema para coste cuadrático desde un punto de vista más geométrico. Estas ideas conducen a una importante cantidad de resultados, por ejemplo, con la interpretación de Otto, algunas de las fórmulas que aparecen en el capítulo \ref{ch:Subcritica}, como \eqref{ec:subcritica_estimacion_a_priori}, resultan mucho más naturales. Para esta Sección seguiremos las notas \cite{VIL_2}.

Para comenzar, en primer lugar presentamos la formulación de Benamou-Brenier del problema de transporte óptimo con coste cuadrático.

Sean $\rho_0$ y $\rho_1$ dos densidades de probabilidad en $\mathbb{R}^n$ con soporte compacto. Pensemos que estamos describiendo dos estados diferentes de una colección de partículas. Para tiempo $t=0$ esta colección de partículas se encuentra en el estado descrito por $\rho_0$, puedes hacer que se muevan alrededor imponiendo cualquier campo vectorial de velocidades en $\mathbb{R}^n$ dependiente del tiempo que desees. Hacer esto conlleva un coste energético que coincide con la energía cinética total de las partículas. Nuestro objetivo es que, en tiempo $t=1$, todas las partículas se encuentren en el estado $\rho_1$, buscando la solución que requiera la menor cantidad de trabajo. Enunciado de manera informal, el objetivo es minimizar la acción
\begin{equation*}
A = \int_0^1 \left( \sum_{x} |\dot{x} (t)|^2 \right) \, dt,
\end{equation*}
con $x$ variando en el conjunto de todas las partículas, distribuido de acuerdo a $\rho_0$ en tiempo $0$, y de acuerdo a $\rho_1$ para tiempo $1$.

Se puede reformular de forma completamente euclídea: Sea $\rho_t$ la densidad de una colección de partículas en tiempo $t$, y $v_t$ su campo de velocidades asociado (definido solo en el soporte de $\rho_t$), el problema anterior ahora consiste en minimizar una energía cinética, 
\begin{equation}\label{ec:Otto_Transporte_optimo}
\inf_{\rho , v} \left\lbrace \int_0^1 \int \rho_t (x) |v_t (x)|^2 \, dx \, dt ; \, \, \, \, \, \, \, \, \, \, \, \, \, \,  \dfrac{\partial \rho_t}{\partial t} + \nabla \cdot (\rho_t v_t) = 0 \right\rbrace , 
\end{equation}
el ínfimo se toma sobre todas las densidades de probabilidad dependientes del tiempo $\left( \rho_t \right)_{0 \leq t \leq 1}$ que coinciden con $\rho_0$ y $\rho_1$ respectivamente para los tiempos $t=0$ y $t=1$, y los campos de velocidades globales dependientes del tiempo $\left( v_t \right)_{0 \leq t \leq 1}$ que dirigen $\rho_t$ expresados según la ecuación de continuidad del lado derecho de \eqref{ec:Otto_Transporte_optimo}.

Según está escrito ahora, este problema de minimización no está definido de forma precisa. Para proporcionarle rigurosidad matemática a la definición uno puede sustituir los valores desconocidos $( \rho , v)$ por $(\rho, P) = (\rho , \rho v)$, con $P$ haciendo referencia a la presión y reescribir \eqref{ec:Otto_Transporte_optimo} como
\begin{equation}\label{ec:Otto_infimo}
\inf_{\rho, P} \int_0^1 \int \dfrac{|P_t (x)|^2}{\rho_t (x)} \, dx \, dt; \, \, \, \, \, \, \, \, \, \, \, \, \, \, \, \, \, \dfrac{\partial \rho_t}{\partial t} + \nabla \cdot P_t =0.
\end{equation}
De esta forma, la minimización del problema tiene sentido bajo ciertas hipótesis de regularidad minimales en $(\rho, P)$ (digamos, por ejemplo, continuidad con respecto al tiempo, $\rho$ tomando valores en el espacio de probabilidad y $P$ en distribuciones vectoriales). Esto viene de juntar la convexidad de $\frac{|P|^2}{\rho}$ como función de $P$ y $\rho$, y la linealidad de las restricciones. Hechas estas anotaciones uno puede escribir el enunciado de manera rigurosa y precisa.

\begin{teor}[Benamou-Brenier, \cite{BRE_2, MCC_2, MCC_3}] \label{th:Benamou_Brenier}
Supongamos que $\rho_0$, $\rho_1$ son densidades de probabilidad en $\mathbb{R}^n$ con soporte compacto, y sea $I$ el valor del ínfimo en \eqref{ec:Otto_infimo}. Entonces,
\begin{equation}\label{ec:Opt_Trans_2-Wasserstein}
I = \min_{T\sharp \rho_0 = \rho_1} \int |x- T(x)|^2 \rho_0 (x) \, dx 
\end{equation}
\end{teor}
La cantidad dada por \eqref{ec:Opt_Trans_2-Wasserstein} se conoce como distancia 2-Wasserstein, que se motiva al formular el Teorema \ref{th:Benamou_Brenier} y que queda definida como
\begin{equation*}
\mathcal{W}_2 (\rho_0, \rho_1 )^2 := \min_{T\sharp \rho_0 = \rho_1} \int |x- T(x)|^2 \rho_0 (x) \, dx .
\end{equation*} 
En estas notas no daremos la demostración del Teorema\footnote{En el Capítulo 8 de \cite{AMB_GIG_SAV} puede encontrarse una demostración del Teorema \ref{th:Benamou_Brenier} de Benamou-Brenier.}, nosotros solo hablaremos de la importancia que este tiene para la teoría de Transporte Óptimo.

Otto descubrió que el Teorema \ref{th:Benamou_Brenier} de Benamou-Brenier llevaba a la re-interpretación de la distancia de Wasserstein como la distancia geodésica para alguna estructura ``Riemanniana'' en el conjunto de densidad de probabilidades:
\begin{equation*}
\mathcal{W}_2 (\rho_0 , \rho_1 )^2 = \inf \left\lbrace \int_0^1 \| \dfrac{\partial \rho}{\partial t} \|^2_{\rho(t)} \, dt; \, \, \, \, \, \, \, \, \, \, \rho (0) = \rho_0, \, \, \rho (1) = \rho_1 , \right\rbrace 
\end{equation*}
donde la ``métrica'' en el ``espacio tangente'' en $\rho$ podría definirse por
\begin{equation*}
\| \dfrac{\partial \rho}{\partial s} \|_{\rho}^2 = \inf_v \left\lbrace \int \rho |v|^2; \, \, \, \, \, \dfrac{\partial \rho}{\partial s} + \nabla \cdot (\rho v) =0 \right\rbrace .
\end{equation*}
La interpretación es: sea $\partial \rho / \partial s$ una variación infinitesimal de la densidad de probabilidad $\rho$. Si pensamos $\rho$ como la probabilidad de un conjunto de partículas, la variación infinitesimal puede asociarse con uno de los muchos campos de velocidades posibles para las partículas, concretamente, con cualquiera de los campos $v$ que satisfacen la ecuación de continuidad $\nabla \cdot ( \rho v) = - \partial \rho / \partial s$. A cada campo de velocidad $v$ se le asocia una energía cinética (la energía cinética total de las partículas), y a la variación infinitesimal le asociamos el mínimo de la energía cinética, de entre todas las posibles elecciones del campo de velocidades. 

Supongamos que $\rho$ es positiva y suave, veamos que un $v$ óptimo debería caracterizarse por ser un \textit{gradiente}. Sea $v$ óptimo y $w$ cualquier campo vectorial con divergencia $0$, entonces $v_{\epsilon} = v + \epsilon (w / \rho)$ es otro campo vectorial admisible, por lo que la energía cinética asociada a $v_{\epsilon}$ no es inferior que la asociada a $v$. De esta forma,
\begin{equation*}
\int \rho |v|^2 \leq \int \rho |v_{\epsilon}|^2 = \int \rho |v|^2 + 2 \epsilon \int v \cdot w + \epsilon^2 \int \dfrac{|w|^2}{\rho} .
\end{equation*}
Simplificando y tomando $\epsilon \rightarrow 0$, vemos que $\int v \cdot w=0$. Dado que $w$ era un vector arbitrario con divergencia $0$, $v$ debería ser un gradiente. 

Obsérvese que, al menos formalmente, la EDP elíptica 
\begin{equation*}
- \nabla \cdot ( \rho \nabla \varphi ) = \dfrac{\partial \rho}{\partial s}
\end{equation*}
debería poder resolverse de forma única (salvo por una constante aditiva) por $\varphi$.

Para resumir, en la formalización de Otto, la distancia de Wasserstein coincide con la distancia geodésica según la métrica definida como: para cada densidad de probabilidad $\rho$, cuando $\partial \rho / \partial s$ y $\partial \rho / \partial t$ sean dos variaciones infinitesimales de $\rho$, entonces 
\begin{equation*}
\left\langle \dfrac{\partial \rho}{\partial s} , \dfrac{\partial \rho}{\partial t} \right\rangle_{\rho} = \int \rho \nabla \varphi \cdot \nabla \psi ,
\end{equation*}
donde $\varphi$, $\psi$ son soluciones respectivamente de
\begin{equation*}
- \nabla \cdot (\rho \nabla \varphi ) = \dfrac{\partial \rho}{\partial s} \, \, \, \, \, \, \, \, \,  \mathrm{y} \, \, \, \, \, \, \, \, - \nabla \cdot ( \rho \nabla \psi ) = \dfrac{\partial \rho}{\partial t}.
\end{equation*}
Además, este formalismo sugiere reglas muy eficientes para realizar los cálculos. De hecho, en cuanto uno ha definido una estructura Riemanniana, la asocia a las reglas del cálculo diferencial que nos permiten profundizar y avanzar en el estudio del flujo gradiente \eqref{ec:flujo_Wasserstein} en la Sección \ref{sec:Sub_JKO}. En particular, se puede definir un gradiente y un operador Hessiano. Recordemos la definición general de un operador gradiente, que se da con la identidad
\begin{equation*}
\left\langle \mathrm{grad}\, E (\rho) , \, \dfrac{\partial \rho}{\partial s} \right\rangle_{\rho} = DE (\rho) \cdot \dfrac{\partial \rho}{\partial s} ,
\end{equation*}
donde $E$ es cualquier función (``suave'') en la variedad que se esté considerando, y $DE$ es el diferencial de $E$. Con las definiciones de arriba, se puede comprobar que el gradiente de las funciones $E$ en el conjunto de las densidades de probabilidad  se define por 
\begin{equation*}
\mathrm{grad}\, E (\rho) = - \nabla \cdot \left( \rho \nabla \dfrac{\delta E}{\delta \rho} \right) .
\end{equation*}
Como consecuencia de esto, ecuaciones con la estructura particular
\begin{equation}\label{ec:Otto_Aplicacion_Importante}
\dfrac{\partial \rho}{\partial t} = \nabla \cdot \left( \rho \nabla \dfrac{\delta E}{\delta \rho} \right)
\end{equation}
pueden verse como flujos gradientes de la energía $E$ con respecto a la estructura diferencial de Otto. Esto es muy interesante si uno mira al comportamiento a largo y corto plazo en ecuaciones como \eqref{ec:Otto_Aplicacion_Importante}. De hecho, existen conexiones fuertes entre el comportamiento del flujo gradiente, $\partial \rho / \partial t = - \mathrm{grad} \, E (\rho)$, y las propiedades de convexidad de $E$.

\chapter{Desigualdades HWI}\label{ap:HWI}

En la sección \ref{sec:Car_McC_Vil_Convergencia_exponencial} hacemos referencia a las desigualdades HWI en repetidas ocasiones para conseguir avanzar con la demostración de convergencia con velocidad exponencial. Por eso les dedicamos este pequeño Apéndice, que nos permite entender un poco más sobre su importancia y su utilidad. Este ha sido redactado siguiendo las notas \cite{MAR_VIL}. 

Estas desigualdades son utilizadas por primera vez por Otto y Villani en \cite{OTT_VIL}. Mezclan la entropía relativa, la distancia de Wasserstein y la información de Fisher relativa, y, a su vez, están relacionadas con el Teorema de Bakry-Emery. El primer caso particular fue probado por Otto en \cite{OTT} para el estudio de ecuaciones del tipo de medio porosos. Podemos dar un enunciado general,
\begin{teor}\label{th:HWI}
Sea $V$ un potencial limitante (confining potential) en $\mathbb{R}^n$, satisfaciendo $D^2 V \geq \lambda I_n $. Asumimos que $e^{-V}$ es una distribución de probabilidad. Entonces, para dos distribuciones de probabilidad cualesquiera $\rho_0$ y $\rho_1$ se tiene la desigualdad HWI
\begin{equation}\label{ec:desig_HWI}
H \left( \rho_0 \, | \, e^{-V} \right) \leq H \left( \rho_1 \, | \, e^{-V} \right) + \mathcal{W}_2 \left( \rho_0 , \rho_1 \right) \sqrt{I \left( \rho_0 \, | \, e^{-V} \right)} - \dfrac{\lambda}{2} \mathcal{W}_2 \left( \rho_0 , \rho_1 \right)^2 .
\end{equation}
\end{teor} 
La desigualdad \eqref{ec:desig_HWI} es bastante fuerte. En particular, si $\lambda > 0$, eligiendo $\rho_0 = e^{-V}$ se obtiene las desigualdades de Talagrand (deducidas en \cite{TAL} por Talagrand desarrollando una idea de Marton \cite{MAR}),
\begin{equation*}
\mathcal{W}_2 \left( \rho_1 , e^{-V} \right) \leq \sqrt{\dfrac{2 H \left( \rho_1 \, | \, e^{-V} \right)}{\lambda}}.
\end{equation*} 
Es decir, el Teorema \ref{th:HWI} permite demostrar que \textit{una desigualdad logarítmica de Sobolev implica una desigualdad de Talgrand}\footnote{Bobkov y Götze fueron los primeros en conjeturar este enunciado \cite{BOB_GOT}. Fue demostrado en \cite{OTT_VIL}.}. Esto es porque el Teorema \ref{th:HWI} incluye en su enunciado las hipótesis necesarias para tener la desigualdad logarítmica de Sobolev. Tomando un caso particular de la conclusión obtenida del Teorema \ref{th:HWI} se recupera la desigualdad de Talagrand.

Por otra parte, la elección de $\rho_1 = e^{-V}$ deja
\begin{equation}\label{ec:des_HWI_particular}
H \left( \rho_0 \, | e^{-V} \right) \leq \mathcal{W}_2 \left( \rho_0, e^{-V} \right) \sqrt{I \left( \rho_0 \, | \, e^{-V} \right)} - \dfrac{\lambda}{2} \mathcal{W}_2 \left( \rho_0, e^{-V} \right)^2 .
\end{equation}
Entonces, si $\lambda > 0$, 
\begin{equation*}
\mathcal{W}_2 \left( \rho_0 , e^{-V} \right) \sqrt{I \left( \rho_0 \, | \, e^{-V} \right)} \leq \dfrac{\lambda}{2} \mathcal{W}_2 \left( \rho_0 , e^{-V} \right)^2 + \dfrac{1}{2 \lambda} I \left( \rho_0 \, | \, e^{-V} \right),
\end{equation*}
por lo que, por \eqref{ec:des_HWI_particular}, implica $H \left( \rho_0 \, | \, e^{-V} \right) \leq \frac{1}{2 \lambda} I \left( \rho_0 \, | \, e^{-V} \right)$, recuperando el Teorema \ref{th:Bakry_Emery} de Bakry-Emery.

\chapter{Ecuación de curvatura Gaussiana} \label{ap:Liouville}
Este Apéndice está motivado por la ecuación de Liouville \eqref{ec:Ecuacion_Liouville} que aparece en el Capítulo \ref{ch:Critica} sobre masa crítica. En estas páginas consideramos dos objetivos principales. Uno de ellos es entender la ecuación para convencernos de la existencia y unicidad de soluciones. El otro es explicar la geometría que rodea a esta ecuación, que a su vez nos permitirá dar las demostraciones de existencia y unicidad para la ecuación de Liouville. Hemos seguido las notas de Chang \cite{CHA} para la elaboración de este texto.

Tenemos la ecuación
\begin{equation}\label{ec:Liouville_Esfera}
\Delta w + e^{2 w} = 1 \, \, \, \mathrm{en} \, \, \mathbb{S}^2 ,
\end{equation}
para la que queremos probar que sus soluciones existen y son únicas. La idea clave del argumento, en la que entra la ecuación de Liouville es que, vamos a interpretar este problema sobre el plano $\mathbb{R}^2$ (así, obtener resultados de existencia y unicidad en una de las dos ecuaciones nos los da para la otra). Utilizamos la proyección estereográfica 
\begin{equation*}
\begin{array}{ll}
\pi : ( \mathbb{S}^n - \mathrm{Polo} \, \, \mathrm{Norte} ) & \rightarrow \mathbb{R}^n \\
\, \, \, \, \, \, \, \, \, \, \, \, \, \, \, \, \, \, \, \, \, \, \, \, \, \, \, \, \, \, \xi & \mapsto x ( \xi)
\end{array}
\end{equation*}
con inversa $\xi = \pi^{-1} (x)$, $\xi_i = \frac{2 x_i}{1 + | x |^2}$, $\xi_{n+1} = \frac{|x|^2 -1}{|x|^2 +1}$, y el cambio de coordenadas,
\begin{equation*}
u (x) = \log \frac{2}{1 + |x|^2} + w ( \xi (x)).
\end{equation*}
Así, para $n=2$, en $\mathbb{R}^2$ la ecuación transformada de \eqref{ec:Liouville_Esfera} pasa a ser
\begin{equation}\label{ec:Liouville_final}
- \Delta u = e^{2u} \, \, \, \mathrm{en} \, \, \mathbb{R}^2.
\end{equation}
Suponiendo $\int_{\mathbb{R}^2} e^{2u} \, dx < \infty$, Chen y Li \cite{CHE_LI} demostraron que \eqref{ec:Liouville_final} es cierto si y solo si $u (x) = \log \frac{2 \lambda}{\lambda^2 + |x-x_0|^2}$, para algún $\lambda > 0$, $x_0 \in \mathbb{R}^2$, utilizando el método de \textit{moving planes}. Esto quiere decir que las bubbles \eqref{ec:Burbujas},
\begin{equation*}
U_{\lambda , \xi} \left( x \right) = \dfrac{1}{\lambda^2} U_0 \left( \dfrac{x - \xi}{\lambda} \right) , \, \, \, U_0 \left( y \right) = \dfrac{8}{( 1 + | y |^2 )^2}, \, \, \, \lambda >0, \, \xi \in \mathbb{R}^2,
\end{equation*}
son las únicas soluciones para la ecuación de Liouville.

También es importante señalar que sin la suposición $\int_{\mathbb{R}^2} e^{2u} \, dx < \infty$, existen otras soluciones analíticas de \eqref{ec:Liouville_final}. De hecho, se tiene una fotografía completa de las ecuaciones en $\mathbb{R}^2$, puede verse la clasificación de Chou y Wan \cite{CHO_WAN}. 

Comprender en mayor profundidad las ideas que hemos expuesto para afirmar que las bubbles son las únicas soluciones de la ecuación de Liouville \eqref{ec:Ecuacion_Liouville} nos anima a estudiar la geometría que envuelve a esta ecuación. 

Para lograr esto hemos divido lo que resta del Apéndice en tres Secciones. En la Sección \ref{ap:sec:Sup_general} atacamos el problema desde un punto de vista más genérico, eligiendo una superficie general y estudiando el problema sobre ella. Tras este análisis previo, en las Secciones \ref{ap:sec:Esfera} y \ref{sec:Esfera_Unicidad} pasamos al caso que más nos interesa, el de la esfera, estudiando el problema en ella, obteniendo existencia y unicidad respectivamente para la solución de la ecuación, terminando de dar forma al enlace que existe con la geometría y comprendiendo por qué obteníamos soluciones con forma de bubbles \eqref{ec:Burbujas}. 

\section{Superficie general}\label{ap:sec:Sup_general}
Sea $(M^2 , g_0)$ una superficie cerrada compacta $2$-dimensional con una métrica dada $g_0$ y curvatura Gaussiana $K_{g_0}$. Estamos interesados en el comportamiento de la curvatura Gaussiana bajo \textit{cambios conformes de la métrica}. Por eso, consideramos la métrica $\bar{g} :=\rho g_0$ para alguna $\rho \in C^{\infty} (M)$, $\rho >0$. Obsérvese que siempre se puede conseguir que $\bar{g}$ y $g_0$ sean coordenadas isotermas, i.e. mientras que la longitud de un vector cambia, el ángulo formado entre cualesquiera dos vectores se preserva bajo el cambio de la métrica de $g_0$ a $\bar{g}$ en $M$. A partir a ahora escribiremos,
\begin{equation*}
\bar{g} = g_w := e^{2 w} g_0
\end{equation*}
para alguna función $w \in C^{\infty} (M)$. Así, tenemos el siguiente resultado que relaciona la curvatura gaussiana de dos métricas conformes:
\begin{prop}
Sea $K_{g_w}$ la curvatura Gaussiana de $(M^2, g_w)$. Entonces
\begin{equation}\label{ec:Ecuacion_Curvatura_Gaussiana}
-\Delta_0 w + K_{g_0} = K_{g_w} e^{2w} .
\end{equation}
\end{prop}
La ecuación \eqref{ec:Ecuacion_Curvatura_Gaussiana} se conoce como la \textit{ecuación de curvatura Gaussiana prescrita}, donde $\Delta_0 = \Delta_{g_0}$ denota el operador de Laplace-Beltrami con respecto a la métrica $g_0$ (a veces denotaremos $\Delta_0$ por $\Delta$).
\begin{proof}
Recordemos la definición del tensor de curvatura Riemanniana. Sea $p \in M^n$, y tomemos una base ortonormal $\left\lbrace e_i \right\rbrace$ del espacio tangente $T_p M$ de $M$ en $p$. Entonces, para dos campos vectoriales $X, Y \in T_p M$ uno tiene
\begin{equation*}
\begin{array}{ll}
Riem(X, Y) & := \nabla_X \nabla_Y - \nabla_Y \nabla_X - \nabla_{[X, Y]} , \\
Riem (e_i , e_j) & \, \, = \nabla_{e_i} \nabla_{e_j} - \nabla_{e_j} \nabla_{e_i} ,
\end{array}
\end{equation*}
donde las dos-formas $Riem$ definen la \textit{curvatura} de la conexión Riemanniana $\nabla$. 

Los \textit{símbolos de Christoffel} de $g$ están dados por 
\begin{equation*}
\Gamma_{ij}^k := \frac{1}{2} g^{kl} \left( \frac{\partial g_{il}}{\partial x^j} + \frac{\partial g_{jl}}{\partial x^i} - \frac{\partial g_{ij}}{\partial x^l} \right),
\end{equation*}
y satisfacen
\begin{equation*}
\nabla_{e_i} e_j = \Gamma_{ij}^k e_k .
\end{equation*}
Sea $R^l_{kij} := g (Riem (e_i, e_j) e_k, e_l)$, entonces el \textit{tensor de Ricci} queda definido como 
\begin{equation*}
Ric_{ij} := Riem^k_{ikj},
\end{equation*}
y la \textit{curvatura escalar} se obtiene otra vez por contracción, quedando
\begin{equation*}
R := Ric_{ij} g^{ij}.
\end{equation*} 
Para $\bar{g} = \rho g_0$, $\rho > 0$, uno calcula directamente (usando $\bar{g}_{il} = \rho (g_0)_{il}$, $\bar{g}^{kl} = \rho^{-1} g_0^{kl}$) que el símbolo de Christoffel $\Gamma_{ij}^k$ de $\bar{g}$ satisface
\begin{equation*}
\bar{\Gamma}_{ij}^k = \Gamma_{ij}^k + \frac{1}{2} \left( \delta_i^k \frac{\partial \log \rho}{\partial x^j} + \delta_j^k \frac{\partial \log \rho}{\partial x^i} - g^{kl} g_{ij} \frac{\partial \log \rho}{\partial x^l} \right).
\end{equation*}
Cuando $n=2$, tomamos $\rho = e^{2 w}$, y después de largos cálculos tediosos que aquí omitimos
\begin{equation*}
\overline{Riem}_{1212} = e^{-2w} ((R_{g_0})_{1212} - 2 \Delta_0 w),
\end{equation*}
que es equivalente a \eqref{ec:Ecuacion_Curvatura_Gaussiana}, ya que $K_{g_0} = \frac{1}{2} (Riem_{g_0})_{1212}$ y $K_{g_w} = \frac{1}{2} \bar{R}_{1212}$. 
\end{proof}

Para terminar de cerrar esta Sección, recurrimos al Teorema de Uniformización, demostrado por Poincaré \cite{POI} y Koebe \cite{KOE_1, KOE_2, KOE_3} y que dice: Toda superficie Riemanniana simplemente conexa es conformemente equivalente al disco unidad abierto, el plano complejo o la esfera de Riemann. En particular, implica que cada superficie de Riemann admite una métrica de Riemann de curvatura constante. Por ello, los casos que más nos interesará estudiar en \eqref{ec:Ecuacion_Curvatura_Gaussiana} son los de curvatura $K_{g_w} = -1$, $0$ y $1$. En lo que queda, vamos a estudiar el problema para curvatura positiva.

\section{Existencia de soluciones en la esfera}\label{ap:sec:Esfera}

Vamos a centrarnos en entender el problema de curvatura prescrita. Centrándonos especialmente en el caso de la esfera con curvatura prescrita $K$. Aquí, Moser dio algunas condiciones de existencia con argumentos variacionales. Consideremos $(M, g) := (\mathbb{S}^2, g_c)$ con la métrica canónica $g_c$ y la curvatura Gaussiana $K_{g_c} \equiv 1$. Entonces, la ecuación de curvatura Gaussiana \eqref{ec:Ecuacion_Curvatura_Gaussiana} se lee como
\begin{equation}\label{ec:Liouville_ap}
\Delta w + K e^{2w} = 1 \, \, \, \, \mathrm{en} \, \, \, \, (\mathbb{S}^2, g_c),
\end{equation}
donde denotamos $\Delta = \Delta_{g_c}$ como antes. Aquí, $K \in C^{\infty}$ es una función dada. En esta Sección comprobaremos existencia para la ecuación.

Moser dio en \cite{MOS} (ver también el libro \cite{MOS_2}) una condición suficiente para resolver \eqref{ec:Liouville_ap}.
\begin{teor}[Moser, \cite{MOS}]\label{th:Moser}
Si $K (- \xi ) =K(\xi )$ para todo $\xi \in \mathbb{S}^2$, y si $\max_{\mathbb{S}^2} K >0$, entonces \eqref{ec:Liouville_ap} tiene una solución $w \in C^{\infty} ( \mathbb{S}^2 )$ con
\begin{equation*}
w (- \xi ) = w (\xi ) \, \, \forall \xi \in \mathbb{S}^2.
\end{equation*} 
\end{teor}
Para la demostración, vamos a considerar una aproximación variacional usando el funcional 
\begin{equation}\label{ec:def:J_K_Moser}
J_K [w] := \log \avint_{\mathbb{S}^2} K e^{2w} \, dv_{g_c} - \frac{1}{4 \pi} \int_{\mathbb{S}^2} | \nabla w |^2 \, dv_{g_c} - 2 \avint_{\mathbb{S}^2} w \, dv_{g_c} ,
\end{equation}
cuyos puntos críticos, $w \in W^{1,2} ( \mathbb{S}^2 )$ satisfacen la ecuación
\begin{equation*}
- \Delta w +1 = \frac{K e^{2w}}{\avint_{\mathbb{S}^2} K e^{2 w} \, dv_{g_c}} \, \, \, \mathrm{en} \, \, \, \mathbb{S}^2 .
\end{equation*}
Entonces, la función desplazada
\begin{equation*}
\tilde{w} := w - \frac{1}{2} \log \avint_{\mathbb{S}^2} K e^{2w} \, dv_{g_c}
\end{equation*}
resuelve \eqref{ec:Liouville_ap}.

Como consecuencia, la demostración se reduce a demostrar la existencia de un punto crítico para el funcional $J_K [ \cdot ]$. Para ello vamos a necesitar los siguientes resultados que daremos sin demostración.
\begin{teor}[Desigualdad de Moser-Trudinger]\label{th:Desigualdad_Moser_Trudinger}
Sea $\Omega \subseteq \mathbb{R}^n$ un dominio acotado, $u \in W_0^{1,n} ( \Omega )$ con $\int_{\Omega} | \nabla u|^n \, dx \leq 1$. Entonces, existe una constante $C = C(n)$, tal que 
\begin{equation*}
\int_{\Omega} e^{\alpha |u|^p} \, dx \leq C |\Omega |,
\end{equation*}
donde $p= \frac{n}{n-1}$, $\alpha \leq \alpha_n := n w_{n-1}^{1/(n-1)}$, $w_k :=$ medida de superficie $k$-dimensional de $\mathbb{S}^k$.
\end{teor}
Para $n=2$ uno tiene $p=2$, $\alpha_2 = 2 w_1 = 4 \pi$. Moser demostró que la constante $\alpha_n$ en el Teorema es sharp. De hecho, construye una sucesión $ u_k \in W_0^{1, n} ( B_1 (0) )$ con $\int_{B_1 (0)} | \nabla u_k |^n \, dx \leq 1$ tal que 
\begin{equation*}
\int_{B_1 (0)} e^{\alpha |u_k|^p} \, dx \rightarrow \infty \, \, \, \, \mathrm{cuando} \, \, \, k \rightarrow \infty,
\end{equation*}
si $\alpha > \alpha_n$.

Damos algún resultado sobre el control con constantes en la métrica $( \mathbb{S}^2 , g_c )$.
\begin{prop}[Moser, \cite{MOS}]\label{prop:Moser_3}
Existe una constante universal $C_1 > 0$, tal que para todo $w \in W^{1,2} ( \mathbb{S}^2 )$ con $\int_{\mathbb{S}^2} | \nabla w|^2 \, d v_{g_c} \leq 1$ y $\int_{\mathbb{S}^2} w \, d v_{g_c} = 0$
\begin{equation}\label{ec:Moser_3}
\int_{\mathbb{S}^2} e^{4 \pi w^2} \, d v_{g_c} \leq C_1.
\end{equation}
\end{prop} 
Uno también puede obtener la cota
\begin{cor}
Para $C_2 := C_1 + \log \frac{1}{4 \pi}$
\begin{equation}\label{ec:Moser_4}
\log \avint_{\mathbb{S}^2} e^{2 w} \, dv_{g_c} \leq \left[ \frac{1}{4 \pi} \int_{\mathbb{S}^2} | \nabla w |^2 \, d v_{g_c} + 2 \avint_{\mathbb{S}^2} w \, d v_{g_c} \right] + C_2
\end{equation}
para todo $w \in W^{1,2} ( \mathbb{S}^2 )$.
\end{cor}
Para $w$ como en la Proposición \ref{prop:Moser_3} con $w$ distinto de $0$, uno puede obtener fácilmente
\begin{equation*}
4 \pi = \int_{\mathbb{S}^2} \, d v_{g_c} < \int_{\mathbb{S}^2} e^{4 \pi w^2} \, d v_{g_c} \leq C_1 ,
\end{equation*}
ya que $C_2 > 0$. Para un dominio en el plano (i.e. $n=2$ en el Teorema \ref{th:Desigualdad_Moser_Trudinger}), Carleson y Chang \cite{CAR_CHA} demostraron la existencia de una función extremal para la desigualdad de Moser-Trudinger para el Teorema \ref{th:Desigualdad_Moser_Trudinger}, observando que la mejor constante $C_2$ en el enunciado del Teorema \ref{th:Desigualdad_Moser_Trudinger} es $> 1 + e$. Este resultado fue extendido por T.L. Soong \cite{SOO} que demostró la existencia de funciones extremales para \eqref{ec:Moser_3} en la Proposición \ref{prop:Moser_3}. 

Para funciones pares en $\mathbb{S}^2$, Moser mejoró sus resultados de la Proposición \ref{prop:Moser_3}:
\begin{prop}[Moser]
Si $w \in W^{1,2} ( \mathbb{S}^2 )$ con $\avint_{\mathbb{S}^2} w \, d v_{g_c} = 0$, $\int_{\mathbb{S}^2} | \nabla w |^2 \, d v_{g_c} \leq 1$ y $w ( \xi ) = w ( - \xi )$ para casi todo $ \xi \in \mathbb{S}^2$, entonces
\begin{equation*}
\int_{\mathbb{S}^2} e^{8 \pi w^2} \, d v_{g_c} \leq C_3.
\end{equation*}
\end{prop}
De nuevo, deducimos
\begin{cor}\label{cor:Moser_5}
Para $C_4 := \log C_3 + \log \frac{1}{4 \pi}$, $a = \frac{1}{2}$,
\begin{equation}\label{ec:Moser_5}
\log \avint_{\mathbb{S}^2} \, d v_{g_c} \leq \left[ a \cdot \frac{1}{4 \pi} \int_{\mathbb{S}^2} | \nabla w |^2 \, d v_{g_c} + 2 \avint_{\mathbb{S}^2} w \, d v_{g_c} \right] + C_4.
\end{equation}
Señalamos que, en lo que respecta a la condición sobre $a$, para las siguientes aplicaciones solo será importante que $a < 1$.
\end{cor}
Con esto podemos volver a la demostración del Teorema \ref{th:Moser}. Dado que $K > 0$ al menos para algún punto, y $K$ es par, 
\begin{equation}\label{ec:Moser_C}
\mathcal{C} := \left\lbrace w \in W^{1,2} ( \mathbb{S}^2) : \int_{\mathbb{S}^2} K e^{2 w} \, d v_{g_c} > 0, \, \, w \, \, \mathrm{par} \, \, \mathrm{a.e.} \right\rbrace \neq \emptyset.
\end{equation}
Consideramos el problema variacional
\begin{equation*}
J_K [ \cdot ] \rightarrow \max J_K [ w ] \, \, \, \mathrm{en} \, \, \mathcal{C},
\end{equation*}
y recordamos que si existe algún $w_0 \in \mathcal{C}$ tal que
\begin{equation*}
\sup_{w \in \mathcal{C}} J_K [w] = J_K [w_0],
\end{equation*}
entonces \eqref{ec:Liouville_ap} tiene solución.

En primer lugar, observamos que $J_K [ \cdot ]$ está acotado por arriba. De hecho, por el Corolario \ref{cor:Moser_5}, \eqref{ec:Moser_5}
\begin{equation*}
\log \avint_{\mathbb{S}^2} K e^{2 w} \, d v_{g_c} \leq \log \max_{\mathbb{S}^2} K + \frac{a}{4 \pi} \int_{\mathbb{S}^2} | \nabla w |^2 \, d v_{g_c} + 2 \avint_{\mathbb{S}^2} w \, d v_{g_c} + C_4,
\end{equation*}
lo que nos lleva a
\begin{equation*}
J_K [ w ] \leq \log \max_{\mathbb{S}^2} K + (a-1) \frac{1}{4 \pi} \int_{\mathbb{S}^2} | \nabla w |^2 \, d v_{g_c} + C_4 < \infty,
\end{equation*}
dado que $a = \frac{1}{2} < 1$. Tomando sucesiones maximizadoras $\left\lbrace w_l \right\rbrace_{l \in \mathbb{N}} \subseteq \mathcal{C}$ con
\begin{equation*}
\lim_{l \to \infty} J_K [ w_l] = \sup_{w \in \mathcal{C}} J_K [w] =: L ,
\end{equation*}
se obtiene,
\begin{equation*}
\begin{array}{ll}
\left( \frac{1-a}{4 \pi} \right) \int_{\mathbb{S}^2} | \nabla w_l |^2 \, d v_{g_c} & \leq \log \max_{\mathbb{S}^2} K + C_4 - J_K [w_l] \\
 & \leq \log \max_{\mathbb{S}^2} K + C_4 + \varepsilon - L
\end{array}
\end{equation*}
para algún $\varepsilon > 0$. Esto implica por la desigualdad de Poincaré que las $w_l$ están uniformemente acotadas en $W^{1,2} ( \mathbb{S}^2 )$, y de ahí, $w_l \rightharpoonup w_0$ en $W^{1,2} (\mathbb{S}^2 )$ para alguna subsucesión. Dado que todas las $w_l$ son pares en casi todo punto, claramente $w_0$ es par para casi todo punto por el Teorema de Rellich\footnote{Se puede leer más sobre la desigualdad de Poincaré y el Teorema de Rellich en la Sección 5.8 del libro \cite{EVA_2}.}. Por la definición de $J_K [ \cdot ]$ en \eqref{ec:def:J_K_Moser}
\begin{equation*}
\left| \log \avint_{\mathbb{S}^2} K e^{2 w_l} \, d v_{g_c} \right| \leq L + C \| w_l \|_{W^{1,2}} \leq \tilde{C} < \infty,
\end{equation*}
y por tanto
\begin{equation*}
\int_{\mathbb{S}^2} K e^{2 ( w_l - \bar{w}_l )} \, d v_{g_c} \geq \min \left\lbrace 4 \pi e^{- \tilde{c} }, 1 \right\rbrace =: c_0 > 0.
\end{equation*}
Usando nuestras convergencias se obtiene
\begin{equation}\label{ec:dem_Moser_final}
\int_{\mathbb{S}^2} K e^{2 w_0} \, d v_{g_c} \geq c_0 > 0,
\end{equation}
y también que
\begin{equation*}
\int_{\mathbb{S}^2} K e^{2 (w_l - \bar{w}_l)} \, dv_{g_c} \to \int_{\mathbb{S}^2} K e^{2(w_0 - \bar{w}_0)} \, dv_{g_c},
\end{equation*}
lo cual implica por \eqref{ec:Moser_C}, que para cualquier $\varepsilon > 0$, existe un $l_0 \in \mathbb{N}$ tal que para todo $l \geq l_0$,
\begin{equation*}
(c_0 - \varepsilon) e^{2 ( \bar{w}_0 - \bar{w}_l )} \leq \int_{\mathbb{S}^2} K e^{2w_0} \, dv_{g_c}.
\end{equation*}
Pero $w_l \to w_0$ en $L^2 (\mathbb{S}^2)$ por el Teorema de Rellich, y de esta forma \eqref{ec:dem_Moser_final} es cierta, con lo que hemos terminado con la demostración del Teorema \ref{th:Moser} y con la parte de existencia de la ecuación \eqref{ec:Liouville_ap}. En la siguiente Sección discutiremos sobre la unicidad.

\section{Unicidad de soluciones en la esfera}\label{sec:Esfera_Unicidad}

Solo queda ver la unicidad para \eqref{ec:Liouville_ap}. Para ello, damos una versión sharp del Corolario \ref{cor:Moser_5}, la desigualdad de Onofri \cite{ONO}.

\begin{teor}[Desigualdad de Onofri]\label{th:Desigualdad_Onofri}
Sea $w \in W^{1,2} (\mathbb{S}^2)$. Entonces
\begin{equation}\label{ec:Onofri_1}
\log \avint_{\mathbb{S}^2} e^{2w} \, dv_{g_c} \leq \frac{1}{4 \pi} \int_{\mathbb{S}^2} | \nabla w |^2 \, dv_{g_c} + 2 \avint_{\mathbb{S}^2} w \, d v_{g_c},
\end{equation}
con igualdad si y solo si,
\begin{equation}\label{ec:Onofri_2}
\Delta w + e^{2 w} = 1,
\end{equation}
i.e.
\begin{equation}\label{ec:Onofri_3}
K_{g_w} \equiv K_{g_c} \equiv 1,
\end{equation}
si y solo si $w = \frac{1}{2} \log |J_{\phi} |$, donde $\phi : \mathbb{S}^2 \rightarrow \mathbb{S}^2$ es una transformación conforme de $\mathbb{S}^2$. En otras palabras, la igualdad en \eqref{ec:Onofri_1} se tiene si y solo si
\begin{equation}\label{ec:Onofri_4}
e^{2w} g_c = \phi^{\ast} (g_c) .
\end{equation}
\end{teor}
La idea clave para demostrar el Teorema \ref{th:Desigualdad_Onofri} es un resultado de Aubin \cite{AUB_2}.
\begin{lema}[Aubin, \cite{AUB_2}]
Sea $\mathcal{S} := \left\lbrace w \in W^{1,2} ( \mathbb{S}^2 ) : \int_{\mathbb{S}^2} e^{2w} x_j \, dv_{g_c} = 0, \, \, j = 1, 2, 3 \right\rbrace $. Entonces para $w \in S$ la siguiente es cierta: Para todo $\varepsilon > 0$ existe una constante $C_{\varepsilon}$ tal que
\begin{equation}\label{ec:lem:Aubin}
\log \avint_{\mathbb{S}^2} e^{2 w} \, dv_{g_c} \leq \left( \frac{1}{2} + \varepsilon \right) \frac{1}{4 \pi} \int_{\mathbb{S}^2} | \nabla w |^2 \, dv_{g_c} + 2 \avint_{\mathbb{S}^2} w \, d v_{g_c} + C_{\varepsilon}.
\end{equation}
\end{lema}
Obsérvese que $\mathcal{S}$ no es muy especial, ya que para cada $w \in C^1 ( \mathbb{S}^2 )$ existe una transformación conforme $\phi : \mathbb{S}^2 \rightarrow \mathbb{S}^2$, tal que
\begin{equation*}
T_{\phi} ( w) := w \, \circ \, \phi + \frac{1}{2} \log | J_{\phi} | \, \, \, \mathrm{en} \, \, \, \mathcal{S}.
\end{equation*}
De hecho, $T_{\phi}$ nos da una correspondencia $1$ a $1$.

Usando \eqref{ec:lem:Aubin}, uno puede obtener compacidad para las sucesiones maximizadoras de $J_K [ \cdot ]$ en $\mathcal{S}$ (ver la definición de $J_K$, \eqref{ec:def:J_K_Moser}). La ecuación de Euler-Lagrange para este problema variacional restringido contiene multiplicadores de Lagrange, que puede demostrarse que pueden hacerse cero usando la condición de Kazdan-Warner \cite{KAZ_WAR}. Finalmente, la unicidad de la solución de \eqref{ec:Onofri_2}, que viene de la ecuación de Euler-Lagrange para $J_K [ \cdot ]$ en $\mathcal{S}$, lleva a $w_{\ast} \equiv 0$ como minimizador. \eqref{ec:Onofri_1} se sigue de $0 = J_K[0] =J_K[w_{\ast}] \leq J_K [w]$ para todo $W^{1,2} ( \mathbb{S}^2 )$, terminando con la demostración.

\chapter{Singularidades tipo I y tipo II}\label{ap:auto-similar}
El estudio de las singularidades es fascinante porque tiene la habilidad de describir una gran variedad de ejemplos provenientes de las ciencias naturales  y otros\footnote{En \cite{KAD}, Kadanoff enumera varias situaciones en las que el estudio de singularidades es relevante.}. El primer ejemplo serían los modelos de quimiotaxia que ya hemos estudiado en el Capítulo \ref{ch:Supercritica}, como \cite{HER_VEL}, que es el caso que nosotros discutimos, o, \cite{BRE_CON_KAD_SCH_VEN}. Otras situaciones en los que estos eventos ocurren son los free-surfaces flows \cite{EGG}, las turbulencias y dinámicas de Euler (singularidades de vórtices tubulares \cite{MOF, GRA_MAR_GER} y laminares \cite{COR_FON_MAN_ROD}), la elasticidad \cite{AUD_BOU}, el condensado Bose-Einstein \cite{BER_RAS}, las ondas físicas no-lineales \cite{MOL_GAE_FIB}, la cosmología de agujeros negros \cite{CHO, MAR_GUN}, o los mercados financieros \cite{SOR}, entre otros.

Para desarrollar el concepto de las singularidades nosotros vamos a seguir las notas de Eggers y Fontelos \cite{EGG_FON}. Empezamos considerando las ecuaciones de evolución
\begin{equation}\label{ec:Ecuacion_de_evolucion}
h_t = F [h],
\end{equation}
donde $F[h]$ representa algún operador diferencial o integral (no-lineal). Para poder centrarnos en los argumentos de la discusión, supondremos que $x$ y $h$ son ambas cantidades escalares, y que las singularidades ocurren en un único punto en tiempo y espacio $x_0$, $T$. Si $\tilde{t} = T - t$ y $\tilde{x} =x-x_0$, buscamos soluciones auto-similares de \eqref{ec:Ecuacion_de_evolucion} que tengan la estructura
\begin{equation}\label{ec:Solucion_ecuacion_evolucion}
h (x, t) = \tilde{t}^{\alpha} H (\tilde{x}/\tilde{t}^{\beta}),
\end{equation}
con los valores adecuadamente elegidos para los exponentes $\alpha$ y $\beta$.

Giga y Kohn \cite{GIG_KOH, GIG_KOH_2} propusieron introducir variables auto-similares $\tau = - \log (\tilde{t} )$ y $y = \tilde{x} / \tilde{t}^{\beta}$ para estudiar el comportamiento asintótico cerca del blow-up. Es decir, poner
\begin{equation}
h (x, t ) = (T-t)^{\alpha} H (y, \tau ),
\end{equation}
y \eqref{ec:Ecuacion_de_evolucion} se convierte en el \textit{sistema dinámico}
\begin{equation}\label{ec:ecuacion_evol_ecuacion_funcional}
H_{\tau} = G [H] := \alpha H - \beta y H_y + F[H].
\end{equation}
Si \eqref{ec:Solucion_ecuacion_evolucion} es efectivamente solución de \eqref{ec:Ecuacion_de_evolucion}, entonces el lado derecho de \eqref{ec:ecuacion_evol_ecuacion_funcional} es independiente de $\tau$ y las soluciones auto-similares de la forma \eqref{ec:Solucion_ecuacion_evolucion} son puntos fijos de \eqref{ec:ecuacion_evol_ecuacion_funcional}. Cuando estudiamos las dinámicas cerca del punto fijo, nos damos cuenta de que utilizando variables auto-similares podemos simplificar el estudio del sistema dinámico \eqref{ec:ecuacion_evol_ecuacion_funcional}, reduciéndolo a un problema de dimensión más baja. La intención que uno tiene al hacer esto, es conseguir información detallada sobre el comportamiento del problema original \eqref{ec:Ecuacion_de_evolucion} cerca del blow-up. 

Con este tipo de estudio se obtiene una forma de \textit{clasificar}, o al menos \textit{caracterizar}, las singularidades. En este Apéndice, las clasificaremos como singularidades tipo I o tipo II. Una clasificación de las singularidades en tipo I y tipo II ha sido expuesta por Barenblatt y Zeldovich en \cite{BAR, BAR_ZEL}. Una solución auto-similar es de tipo I si \eqref{ec:Solucion_ecuacion_evolucion} solo resuelve \eqref{ec:Ecuacion_de_evolucion} para un conjunto de exponentes $\alpha$, $\beta$; sus valores pueden fijarse, por ejemplo, por medio de un análisis dimensional o por simetría (y por consiguiente, son radiales). Las soluciones son de tipo II si para \eqref{ec:Solucion_ecuacion_evolucion} existen soluciones localmente para un conjunto continuo de exponentes $\alpha$, $\beta$. No obstante, en general, estas soluciones son inconsistentes con las condiciones de frontera o las condiciones iniciales. Imponerlas nos lleva a un problema no-lineal de autovalores, cuyas soluciones toman exponentes irracionales en general.

Las técnicas involucradas en este tipo de argumentos son bastante complejas, así que nosotros aspiramos a explicarlas a través de ejemplos. En la Sección \ref{ap:sec:Tipo_I}, analizamos los blow-up tipo I estudiando la difusión superficial de átomos. Tras esto, en la Sección \ref{ap:sec:Tipo_II}, intentamos comprender las singularidades tipo II a través de la ecuación de Burgers.

\section{Blow-up tipo I}\label{ap:sec:Tipo_I}
El ejemplo que presentamos sobre singularidades tipo I se debe a Barenblatt \cite{BAR}. Este es el de una superficie sólida evolucionando bajo la acción de difusión superficial en la que los átomos migran a lo largo de la superficie, conducidos por gradientes de potenciales químicos. Las ecuaciones para el caso de simetría con respecto a un eje (donde la superficie libre es descrita por el radio $h(x, t)$), son \cite{NIC_MUL}
\begin{equation}\label{ec:Ecuacion_Atomos}
h_t = \frac{1}{h} \left[ \frac{h}{(1 +h_x^2)^{1/2}} \kappa_x \right],
\end{equation}
donde,
\begin{equation}\label{ec:Tipo_I_Curvatura_media}
\kappa = \frac{1}{h (1 + h_x^2 )^{1/2}} - \frac{h_{xx}}{(1+h_x^2)^{3/2}}
\end{equation}
es la curvatura media. En \eqref{ec:Ecuacion_Atomos}, \eqref{ec:Tipo_I_Curvatura_media}, todas las longitudes se han construido de forma adimensional, usando una escala exterior $R$ (como por ejemplo, el radio inicial), y la escala de tiempo $R^4/D_4$, donde $D_4$ es una constante de difusión de cuarto orden. 

Para tiempo $\tilde{t} \ll 1$, lejos de la desintegración, un análisis dimensional implica que $l = \tilde{t}^{1/4}$ es una escala de longitud local. Esto sugiere la forma auto-similar
\begin{equation}\label{ec:Atomos_similitud}
h (x, t) = \tilde{t}^{1/4} H ( \tilde{x} / \tilde{t}^{1/4} ),
\end{equation}
y así, los exponentes $\alpha$, $\beta$ de \eqref{ec:Solucion_ecuacion_evolucion} son fijados a través de un análisis dimensional, lo cuál es típico para auto-similitud de tipo I. Así, la forma similar de la EDP resulta
\begin{equation}\label{ec:Atomos_forma_similar}
- \frac{1}{4} (H - y H_y ) = \frac{1}{H} \left[ \frac{H}{(1+ H_y^2)^{1/2}} \kappa_y \right]_y , \, \, \, \, y = \frac{\tilde{x}}{\tilde{t}^{1/4}}
\end{equation}
donde $\kappa$ es la curvatura de $H$. 

Las soluciones de \eqref{ec:Atomos_forma_similar} han sido estudiadas extensivamente en \cite{BER_BER_WIT}. Queremos utilizar técnicas de matching, que, a grandes rasgos, consisten en estudiar por separado el problema cerca de la singularidad y lejos de esta y, tras ello, unir ambas partes. Para asegurar el matching a una solución exterior independiente del tiempo, la dependencia del término principal en tiempo debe salir de \eqref{ec:Atomos_similitud}, implicando que 
\begin{equation}\label{ec:Atomos_cond_creciemiento}
H (y) \sim c | y |, \, \, \, \, y \rightarrow \pm \infty.
\end{equation}
En general, la condición para obtener este matching para soluciones auto-similares de la forma \eqref{ec:Solucion_ecuacion_evolucion} es
\begin{equation}\label{ec:condicion_matching}
H (y) \sim c |y|^{\frac{\alpha}{\beta}}, \, \, \, \, y \rightarrow \pm \infty.
\end{equation}
Todas las soluciones de la ecuación de similitud \eqref{ec:Ecuacion_Atomos}, y que siguen una condición de crecimiento \eqref{ec:Atomos_cond_creciemiento} son simétricas, y forman un conjunto infinito discreto. El hecho de que las soluciones similares permitidas formen un conjunto discreto implica un gran acuerdo de ``universalidad'' en la forma en la que los pinching pueden ocurrir. Esto significa que la solución local es independiente de la solución exterior, y además que la primera impone restricciones sobre la segunda. En particular, el prefactor $c$ en \eqref{ec:Atomos_cond_creciemiento} debe determinarse como parte de la solución.

Otro ejemplo interesante de soluciones que presentan estas singularidades tipo I son los \textit{thin films and thin jets}, estudiadas por Chou y Kwong en \cite{CHO_KWO}, y que son una generalización de la \textit{long-wave thin-film equation}
\begin{equation}\label{ec:long-wave_thin-film}
h_t + \nabla \cdot (h^n \nabla \Delta h - B h^m \nabla h ) =0, \, \, \, \, n > 0.
\end{equation}
El flujo de masa en esta ecuación tiene dos contribuciones. La primera se debe a la tensión de la superficie, y la segunda es a causa de un potencial externo. Cuando $n = m = 3$, entonces $z = h(x, t)$ representa la altura de una capa o una gota de un fluido viscoso sobre una superficie plana, localizada en $z=0$ y donde el potencial externo es la gravedad. Si $B$ es negativo, \eqref{ec:long-wave_thin-film} describe una capa que cuelga de un techo. Independientemente del signo de $B$, en esta situación, \eqref{ec:long-wave_thin-film} no tiene singularidad.

Decimos que las soluciones de \eqref{ec:long-wave_thin-film} desarrollan puntos singulares si $h$ va a cero en tiempo finito. Esto ocurre si uno añade fuerzas de van der Waals, que en el término principal implican $n=3$ y $m=-1$ con $B < 0$. Por medio de métodos numéricos se muestra para este caso la existencia de soluciones auto-similares que tocan el suelo con simetría radial de la forma
\begin{equation*}
h(r, t) = \tilde{t}^{\frac{1}{5}} H( y), \, \, \, y = r/ \tilde{t}^{\frac{2}{5}}.
\end{equation*}
Soluciones auto-similares que tocan a lo largo de una recta también existen, pero son inestables.

\section{Blow-up tipo II}\label{ap:sec:Tipo_II}
En el ejemplo de la Sección anterior, los exponentes pueden determinarse recurriendo a un análisis dimensional, o gracias a argumentos de simetría, y, por tanto, se pueden asumir valores racionales. No obstante, en otros problemas el comportamiento del reescalamiento depende de parámetros externos, fijados, por ejemplo, por las condiciones iniciales. En este caso, en el argumento de reescalamiento podemos asumir cualquier valor. Normalmente, este queda fijado por una condición de compatibilidad, que al final resulta en una respuesta irracional. Llamaremos a esta situación auto-similaridad tipo II. Dado que es relativamente raro poder seguir los resultados de forma analítica, en esta Sección, tal y como señalábamos al comienzo del Apéndice, vamos a centrarnos en explicar el blow-up tipo II usando un ejemplo.

Consideramos el modelo más sencillo de formación de una ola de choque en dinámica de gases, que es la ecuación de Burgers
\begin{equation}\label{ec:Ecuacion_Burgers}
u_t + u_{xx} = 0.
\end{equation}
Normalmente, se cree que cualquier sistema de leyes de conservación que exhibe blow-up se comportará localmente como \eqref{ec:Ecuacion_Burgers} \cite{ALI}. 

Durante esta revisión, solo consideramos las dinámicas hasta la singularidad. Esta estructura que surge después de la singularidad depende de la regularización usada, ya que la continuación en tiempos posteriores a la singularidad no es única \cite{DAF, WHI}. 

Se sabe \cite{LAN_LIF} que \eqref{ec:Ecuacion_Burgers} puede resolverse de forma exacta usando el método de las características. Este método consiste en notar que la velocidad permanece constante a lo largo de la curva característica
\begin{equation*}
z = u_0 (x) t + x,
\end{equation*}
donde $u_0 (x) = u (x, 0)$ es la condición inicial. De esta forma,
\begin{equation}\label{ec:Burgers_condicion_inicial}
u (z, t) = u_0 (x)
\end{equation}
es una solución exacta de \eqref{ec:Ecuacion_Burgers}, dada implícitamente.

Es geometricamente obvio que cuando $u_0 (x)$ tiene una pendiente negativa, las características se cruzarán en tiempo finito y producirán una discontinuidad de la solución. Esto ocurre cuando $\partial z / \partial x = 0$, lo cuál sucede por primera vez en el tiempo de singularidad
\begin{equation*}
T = \min \left\lbrace - \frac{1}{\partial_x u_0 (x)} \right\rbrace ,
\end{equation*}
en una posición espacial $x = x_m$. Esto significa que una singularidad se formará por primera vez en
\begin{equation*}
x_0 = x_m - \frac{u_0 (x_m)}{\partial_x u_0 (x_m)} .
\end{equation*}
Dado que \eqref{ec:Ecuacion_Burgers} es invariante bajo cualquier cambio de velocidad, podemos asumir sin pérdida de generalidad que $u_0 (x_m) = 0$, y de forma que $x_0 = x_m$. Esto significa que la velocidad es cero en la singularidad. Lo siguiente, es analizar su formación usando la coordenada local $\tilde{x}$, $\tilde{t}$. En \cite{POM_LEB_GUY_GRI}, Pomeau, Le Berre, Guyenne y Grilli hicieron esto expandiendo la condición inicial \eqref{ec:Burgers_condicion_inicial} usando ideas de la teoría de las catástrofes \cite{ARN}. Nosotros, en cambio, recurriremos a ideas de auto-similaridad para poder entender estas técnicas.

El comportamiento local de \eqref{ec:Ecuacion_Burgers} cerca de $T$ puede obtenerse por medio de una reescalamiento
\begin{equation*}
u (x, t) = \tilde{t}^{\alpha} U ( \tilde{x} / \tilde{t}^{\alpha + 1} ),
\end{equation*}
que resuelve \eqref{ec:Ecuacion_Burgers}, dando una solución para la ecuación de auto-similaridad
\begin{equation*}
- \alpha U + (1 + \alpha ) y U_y + U U_y = 0,
\end{equation*}
con solución implícita
\begin{equation}\label{ec:Burgers_solucion_implicita}
y= - U - C U^{1+1/\alpha}.
\end{equation}
El caso especial $\alpha = 0$ tiene solución $U = -y$, que es inconsistente con la condición de matching \eqref{ec:condicion_matching}, y por ello, esta es descartada.

Nos quedamos con un continuo de posibles exponentes de reescalamiento $\alpha > 0$, algo que es típico para las auto-similitudes de tipo II. No obstante, una sucesión discreta infinita de exponentes $\alpha_n$ que es elegida según los requisitos impuestos por \eqref{ec:Burgers_solucion_implicita} define una función suave para todo $y$. Concretamente, uno tiene $1 + 1/\alpha$ impar, o
\begin{equation*}
\alpha_i = \frac{1}{2 i + 2} , \, \, \, \, i = 0, 1, 2, \cdots ,
\end{equation*}
y denotamos el correspondiente perfil de auto-similaridad para $U_i$. La constante $C$ en \eqref{ec:Burgers_solucion_implicita} debe ser positiva pero arbitraria. Queda fijada por las condiciones iniciales, otro sello distintivo de las auto-similitudes tipo II. No obstante, $C$ puede normalizarse a $1$ al reescalar $x$ y $U$. La solución con $\alpha_0$,
\begin{equation}\label{ec:Burgers_singularidad}
u (x, t) \tilde{t}^{1/2} U_0 ( \tilde{x} / \tilde{t}^{3/2} ),
\end{equation} 
es la única que es estable, todas las soluciones de orden más alto son inestables.

Sigue siendo interesante mirar a algunas posibles excepciones a la forma de blow-up dada arriba,
\begin{equation}\label{ec:Burgers_modificado}
u_t + u u_x = u^{\sigma }
\end{equation}
Esta ecuación también se resuelve fácilmente usando características. Para $\sigma \leq 2$ se esperan dos comportamientos posibles. Para pequeños datos iniciales $u_0 (x)$ todavía se forma una singularidad como \eqref{ec:Burgers_singularidad}, aunque hay que resaltar que puede que $u$ también se vaya a infinito. No obstante, existe una frontera entre los dos comportamientos \cite{ALI}, donde la pendiente explota al mismo tiempo que $u$ va a infinito. En este caso, uno espera que todos los términos en \eqref{ec:Burgers_modificado} tengan el mismo orden, dando
\begin{equation}\label{ec:Burgers_modificado_bolw-up}
u (x, t) = \tilde{t}^{\frac{1}{1- \sigma}} U (y), \, \, \, \, y = \tilde{x} / \tilde{t}^{\frac{\sigma - 2}{\sigma - 1}},
\end{equation}
con ecuación de auto-similitud
\begin{equation}\label{ec:Burgers_modificado_auto-similar}
\frac{U}{1 - \sigma} + \frac{\sigma - 2}{\sigma - 1} y U_y = U^{\sigma} - U U_y.
\end{equation}
La ecuación de \eqref{ec:Burgers_modificado_auto-similar} que tiene el decaimiento correcto en el infinito es
\begin{equation}\label{ec:Solucion_Burgers_modificado}
y = - \frac{1}{(\sigma -2) U^{\sigma -2}} \pm C \frac{(1- (\sigma -1) U^{\sigma -1})^{\frac{\sigma -2}{\sigma -1}}}{U^{\sigma -2}},
\end{equation}
donde $C>0$ es una constante arbitraria. Los signos $+$ y $-$ describen la solución a la derecha y la izquierda de $y^{\ast} = - (\sigma -1)^{\frac{\sigma -2}{\sigma -1}} / (\sigma -2)$, respectivamente. La solución auto-similar \eqref{ec:Solucion_Burgers_modificado} no es suave en su máximo, su primera derivada se comporta como $U_y \varpropto (y - y^{\ast})^{1/(\sigma -2)}$. Esto puede entenderse de la solución exacta, para que el blow-up ocurra en el momento en el que un shock es formado, el perfil es inicial debe tener un máximo con la misma regularidad que \eqref{ec:Solucion_Burgers_modificado}. Así pues, la situación dirigida a \eqref{ec:Burgers_modificado_bolw-up} es una muy especial, que requiere condiciones iniciales muy peculiares.

\chapter*{Agradecimientos}
Estas notas constituyen el trabajo de fin de grado presentado por el autor en junio del año 2021 para finalizar el grado en matemáticas ofertado por la Universidad Autónoma de Madrid. Por ello, como autor, quiero agradecer a la profesora María del Mar González, la tutora de este trabajo de fin de grado, por su paciencia al enseñarme todo lo que sé sobre las ecuaciones de Keller-Segel, por transmitirme su pasión por las matemáticas y por descubrirme lo maravillosa que pueden ser, con charlas interminables sobre EDPs y... sobre la vida en general. Además, el autor es apoyado por la ``Advanced Grant Nonlocal-CPD (Nonlocal PDEs for Complex Particle Dynamics: Phase Transitions, Patterns and Synchronization)'' de la European Research Council Executive Agency (ERC) bajo el programa de investigación e innovación European Union's Horizon 2020 (acuerdo de beca No. 883363).

\cleardoublepage

\end{document}